\documentclass[11pt]{article}
\usepackage[mathscr]{eucal}
\usepackage{amsmath}

\newcommand{\ds}{\displaystyle}
\DeclareMathAlphabet{\mathpzc}{OT1}{pzc}{m}{it}

\newcommand{{\M}}{{\textsf{{L}$_\sharp^{2}$}}}

\newcommand{\bsfH}{\boldsymbol{\mathsf H}}

\newcommand{\bsfU}{\boldsymbol{\mathsf U}}

\usepackage{bbm}
\usepackage{dsfont}
\usepackage{color}
\usepackage{mathrsfs}
\usepackage{amssymb}
\usepackage{calligra}
\usepackage{fontenc}
\usepackage{amsbsy}
\usepackage{graphicx}
\graphicspath{%
    {converted_graphics/}
    {/}
}
\begin{document}

\newcommand{\sfU}{{\sf U}}
\newcommand{\sfF}{{\sf F}}
\newcommand{\sfW}{{\sf W}}
\newcommand{\sfD}{{\sf D}}
\newcommand{\sfL}{{\sf L}}
\newcommand{\sfu}{{\sf u}}
\newcommand{\sfw}{{\sf w}}
\newcommand{\sfz}{{\sf z}}
\newcommand{\sfp}{{\sf p}}
\newcommand{\sfv}{{\sf v}}
\newcommand{\simge}{\ba{cc}\vspace*{-2.4mm}>\\ \sim\ea }
\newcommand{\simle}{\ba{cc}\vspace*{-2.4mm}<\\ \sim\ea }
\newcommand{\Cdot}{\!\cdot\!}
\newcommand{\sq}{{$\sqcap\!\!\!\!\sqcup$}}
\newcommand{\Eu}{{\rm I\,\!\! E}}
\newcommand{\Io}{\Int{\Omega}{}}
\newcommand{\Id}{\Int{\cald}{}}
\newcommand{\Div}{\mbox{\rm div}\,}
\newcommand{\tr}{\mbox{\rm tr}\,}
\newcommand{\grad}{\mbox{\rm grad}\,}
\newcommand{\supp}{\mbox{\rm supp}\,}
\newcommand{\curl}{\mbox{\rm curl}\,}
\newcommand{\Ido}{\Int{\partial\Omega}{}}
\newcommand{\IdS}{\Int{\Sigma}{}}
\newcommand{\Oint}[2]{{\displaystyle \oint_{#1}^{#2}}}
\newcommand{\Int}[2]{{\displaystyle \int_{ #1}^{ #2}}}
\newcommand{\Lim}[1]{{\displaystyle \lim_{ #1}}}
\newcommand{\Limsup}[1]{{\displaystyle \limsup_{\footnotesize #1}}}
\newcommand{\Liminf}[1]{{\displaystyle \liminf_{\footnotesize #1}}}
\newcommand{\Sup}[1]{{\displaystyle \sup_{#1}}}
\newcommand{\Inf}[1]{{\displaystyle \inf_{#1}}}
\newcommand{\Max}[1]{{\displaystyle \max_{#1}}}
\newcommand{\Min}[1]{{\displaystyle \min_{#1}}}
\newcommand{\Sum}[2]{{\displaystyle \sum_{#1}^{#2}}}
\newcommand{\Prod}[2]{{\displaystyle \prod_{#1}^{#2}}}
\newcommand{\BCup}[2]{{\displaystyle \bigcup_{#1}^{#2}}}
\newcommand{\BCap}[2]{{\displaystyle \bigcap_{#1}^{#2}}}
\newcommand{\Frac}[2]{\displaystyle{\frac{\displaystyle{#1}}{\displaystyle{#2}}}}
\newcommand{\norm}[1]{\left\|{#1}\right\|}
\newcommand{\Norm}[1]{\langle\langle{#1}\rangle\rangle_q}
\newcommand{\No}[1]{\langle\!\langle{#1}\rangle\!\rangle}
\newcommand{\NO}[1]{{\langle{#1}\rangle}_{\lambda,q}}
\newcommand{\beea}{\begin{eqnarray}}
\newcommand{\eeea}{\end{eqnarray}}
\newcommand{\ms}{\medskip\smallskip}
\newcommand{\bs}{\bigskip}
\newcommand{\ps}{\par\smallskip}
\newcommand{\bfe}{{\mbox{\boldmath $e$}} }
\newcommand{\pni}{{\par\noindent}}
\newcommand{\bfq}{{\mbox{\boldmath $q$}} }
\newcommand{\bfz}{{\mbox{\boldmath $z$}} }
\newcommand{\0}{{\mbox{\boldmath $0$}} }
\newcommand{\LE}{\!\!\!&\le&\!\!\!}
\newcommand{\BL}[1]{{\par\smallskip{\bf Lemma #1.}}}
\newcommand{\BT}[1]{{\par\smallskip{\bf Theorem #1.}}}
\newcommand{\Ln}{[\!|}
\newcommand{\Rn}{|\!]}
\newcommand{\n}[1]{{\Ln{#1}\Rn}} 
\newcommand{\nq}[1]{{\Ln{#1}\Rn}_{q}} 
\newcommand{\nqr}[1]{{\Ln{#1}\Rn}_{q,r}} 
\newcommand{\Nq}[1]{{\langle{#1}\rangle}_{q}} 
\newcommand{\Nql}[1]{{\langle{#1}\rangle}_{\lambda,q}} 
\newcommand{\Nqr}[1]{{\langle{#1}\rangle}_{q,r}}
\newcommand{\N}[1]{{|\!\!|\!\!|\,{#1}\,|\!|\!\!|_2}}
\newcommand{\EA}[2]{$$#1$$%
\vspace{-6.mm}
\begin{equation}
\end{equation}
\vspace{-6.mm}
$$
#2
\setlength{\belowdisplayskip}{3mm}
\setlength{\belowdisplayshortskip}{3mm}
$$
}
\newcommand{\A}[2]{$$#1$$%
\vspace{-4.mm}
$$
#2
\setlength{\belowdisplayskip}{3mm}
\setlength{\belowdisplayshortskip}{3mm}
$$
}
\newcommand{\BF}{\begin{footnotesize}}
\newcommand{\EF}{\end{footnotesize}}
\setlength{\jot}{.15in}
\newcommand{\pde}[2]{{\displaystyle \frac{\mbox{$\partial #1$}}{\mbox{$\partial #2$}}}}
\newcommand{\ode}[2]{{\displaystyle \frac{\mbox{$d #1$}}{\mbox{$d #2$}}}}
\newcommand{\f}[2]{\frac{\mbox{$#1$}}{\mbox{$ #2$}}}
\newcommand{\bi}{\begin{itemize}}
\newcommand{\ei}{\end{itemize}}
\newcommand{\ed}{\end{document}}
\newcommand{\be}{\begin{equation}}
\newcommand{\ba}{\begin{array}}
\newcommand{\ea}{\end{array}}
\newcommand{\ee}{\end{equation}}
\newcommand{\eeq}[1]{\label{eq:#1}\end{equation}}
\newcommand{\real}{{\mathbb R}}
\newcommand{\compl}{{\mathbb C}}
\def\Id{\mbox{\boldmath $1$}}
\def\zero{\mbox{\boldmath $0$}}
\newcommand{\PP}{{\rm I\!\!\,P}}
\newcommand{\nat}{{\mathbb N}}
\newcommand{\bfpsi}{\mbox{\boldmath $\psi$}}
\newcommand{\bfomega}{\mbox{\boldmath $\omega$}}
\newcommand{\bfvaromega}{\mbox{\boldmath $\varpi$}}
\newcommand{\bfOmega}{\mbox{\boldmath $\Omega$}}
\newcommand{\bfTheta}{\mbox{\boldmath $\Theta$}}
\newcommand{\bfmu}{\mbox{\boldmath $\mu$}}
\newcommand{\bfx}{\mbox{\boldmath $x$}}
\newcommand{\bfy}{\mbox{\boldmath $y$}}
\newcommand{\bfPsi}{\mbox{\boldmath $\Psi$}}
\newcommand{\bfxi}{\mbox{\boldmath $\xi$}}

\newcommand{\bfphi}{\mbox{\boldmath $\varphi$}}
\newcommand{\bfhi}{\mbox{\boldmath $\phi$}}
\newcommand{\bfPhi}{\mbox{\boldmath $\Phi$}}
\newcommand{\bfv}{{\mbox{\boldmath $v$}} }
\newcommand{\bfu}{{\mbox{\boldmath $u$}} }
\newcommand{\bfuf}{{\mbox{\footnotesize\boldmath $u$}} }
\newcommand{\bfw}{{\mbox{\boldmath $w$}} }
\newcommand{\bff}{{\mbox{\boldmath $f$}} }
\newcommand{\bfa}{{\mbox{\boldmath $a$}} }
\newcommand{\bfi}{{\mbox{\boldmath $i$}} }
\newcommand{\bfj}{{\mbox{\boldmath $j$}} }
\newcommand{\bfc}{{\mbox{\boldmath $c$}} }
\newcommand{\bfo}{{\mbox{\boldmath $o$}} }
\newcommand{\bfp}{{\mbox{\boldmath $p$}} }
\newcommand{\bfkp}{{\mbox{\footnotesize{\boldmath $k$}}} }
\newcommand{\bfka}{{\mbox{\footnotesize{\boldmath $k^*$}}} }
\newcommand{\bft}{{\mbox{\boldmath $t$}} }
\newcommand{\bfd}{{\mbox{\boldmath $d$}} }
\newcommand{\bfl}{{\mbox{\boldmath $l$}} }
\newcommand{\bfr}{{\mbox{\boldmath $r$}} }
\newcommand{\bfk}{{\mbox{\boldmath $k$}} }
\newcommand{\bfA}{{\mbox{\boldmath $A$}} }
\newcommand{\bfS}{{\mbox{\boldmath $S$}} }
\newcommand{\bfO}{{\mbox{\boldmath $O$}} }
\newcommand{\bfM}{{\mbox{\boldmath $M$}} }
\newcommand{\bfP}{{\mbox{\boldmath $P$}} }
\newcommand{\bfB}{{\mbox{\boldmath $B$}} }
\newcommand{\bfR}{{\mbox{\boldmath $R$}} }
\newcommand{\bfC}{{\mbox{\boldmath $C$}} }
\newcommand{\bfD}{{\mbox{\boldmath $D$}} }
\newcommand{\bfQ}{{\mbox{\boldmath $Q$}} }
\newcommand{\bfZ}{{\mbox{\boldmath $Z$}} }
\newcommand{\bfG}{{\mbox{\boldmath $G$}} }
\newcommand{\bfE}{{\mbox{\boldmath $E$}} }
\newcommand{\bfX}{{\mbox{\boldmath $X$}} }
\newcommand{\bfY}{{\mbox{\boldmath $Y$}} }
\newcommand{\bfH}{{\mbox{\boldmath $H$}} }
\newcommand{\bfI}{{\mbox{\boldmath $I$}} }
\newcommand{\bfJ}{{\mbox{\boldmath $J$}} }
\newcommand{\bfN}{{\mbox{\boldmath $N$}} }
\newcommand{\bfh}{{\mbox{\boldmath $h$}} }
\newcommand{\bfm}{{\mbox{\boldmath $m$}} }
\newcommand{\bfone}{{\mbox{\boldmath $1$}} }
\newcommand{\hs}{{\rm I}\!\!\,{\rm R}^3_+}
\newcommand{\cala}{{\cal A}}
\newcommand{\calb}{{\cal B}}
\newcommand{\calc}{{\cal C}}
\newcommand{\cald}{{\cal D}}
\newcommand{\cale}{{\cal E}}
\newcommand{\calf}{{\cal F}}
\newcommand{\calg}{{\cal G}}
\newcommand{\calh}{{\cal H}}
\newcommand{\cali}{{\cal I}}
\newcommand{\calj}{{\cal J}}
\newcommand{\calk}{{\cal K}}
\newcommand{\call}{{\cal L}}
\newcommand{\calm}{{\cal M}}
\newcommand{\caln}{{\cal N}}
\newcommand{\calo}{{\cal O}}
\newcommand{\calp}{{\cal P}}
\newcommand{\calq}{{\cal Q}}
\newcommand{\calr}{{\cal R}}
\newcommand{\cals}{{\cal S}}
\newcommand{\calt}{{\cal T}}
\newcommand{\calu}{{\cal U}}
\newcommand{\calv}{{\cal V}}
\newcommand{\calx}{{\cal X}}
\newcommand{\caly}{{\cal Y}}
\newcommand{\calw}{{\cal W}}
\newcommand{\calz}{{\cal Z}}
\newcommand{\bfsigma}{\mbox{\boldmath $\sigma$}}
\newcommand{\bfSigma}{\mbox{\boldmath $\Sigma$}}
\newcommand{\bftau}{\mbox{\boldmath $\tau$}}
\newcommand{\bfeta}{\mbox{\boldmath $\eta$}}
\newcommand{\bfT}{{\mbox{\boldmath $T$}} }
\newcommand{\bfV}{{\mbox{\boldmath $V$}} }
\newcommand{\bfU}{{\mbox{\boldmath $U$}} }
\newcommand{\bfW}{{\mbox{\boldmath $W$}} }
\newcommand{\bfF}{{\mbox{\boldmath $F$}} }
\newcommand{\bfK}{{\mbox{\boldmath $K$}} }
\newcommand{\bfL}{{\mbox{\boldmath $L$}} }
\newcommand{\bfb}{{\mbox{\boldmath $b$}} }
\newcommand{\bfg}{{\mbox{\boldmath $g$}} }
\newcommand{\bfn}{{\mbox{\boldmath $n$}} }
\newcommand{\bfs}{{\mbox{\boldmath $s$}} }
\newcommand{\cf}{{\it cf.} }
\newcommand{\io}{\int_\Omega}
\newcommand{\1}{\item[({\it i})]}
\newcommand{\2}{\item[({\it ii})]}
\newcommand{\3}{\item[({\it iii})]}
\newcommand{\4}{\item[({\it iv})]}
\newcommand{\5}{\item[({\it v})]}
\newcommand{\6}{\item[({\it vi})]}
\newcommand{\7}{\item[({\it vii})]}
\newcommand{\8}{\item[({\it viii})]}
\newcommand{\9}{\item[({\it xi})]}
\newcommand{\ido}{\int_{\partial\Omega}}
\newcommand{\half}{\mbox{$\frac{1}{2}$}}
\def\parallel{\|}
\def\mid{|}
\def\Bbb R{\real}
\def\hat{\widehat}
\def\tilde{\widetilde}
\def\bar{\overline}
\newcommand{\threehalves}{3\over 2}
\newcommand{\bfPi}{\mbox{\boldmath $\Pi$}}
\newcommand{\bfchi}{\mbox{\boldmath $\chi$}}
\newcommand{\bfalpha}{\mbox{\boldmath $\alpha$}}
\newcommand{\bfbeta}{\mbox{\boldmath $\beta$}}
\newcommand{\bfgamma}{\mbox{\boldmath $\gamma$}}
\newcommand{\bfdelta}{\mbox{\boldmath $\delta$}}
\newcommand{\bfzeta}{\mbox{\boldmath $\zeta$}}
\newcommand{\bfUpsilon}{\mbox{\boldmath $\Upsilon$}}
\newcommand{\bfGamma}{\mbox{\boldmath $\Gamma$}}
\newcommand{\bfcala}{\mbox{\boldmath ${\cal A}$}}
\newcommand{\bfcalm}{\mbox{\boldmath ${\cal M}$}}
\newcommand{\bfcaln}{\mbox{\boldmath ${\cal N}$}}
\newcommand{\bfcalq}{\mbox{\boldmath ${\cal Q}$}}
\newcommand{\bfcalb}{\mbox{\boldmath ${\cal B}$}}
\newcommand{\bfcalc}{\mbox{\boldmath ${\cal C}$}}
\newcommand{\bfcali}{\mbox{\boldmath ${\cal I}$}}
\newcommand{\bfcalg}{\mbox{\boldmath ${\cal G}$}}
\newcommand{\bfcalh}{\mbox{\boldmath ${\cal H}$}}
\newcommand{\bfcalk}{\mbox{\boldmath ${\cal K}$}}
\newcommand{\bfcalt}{\mbox{\boldmath ${\cal T}$}}
\newcommand{\bfcalx}{\mbox{\boldmath ${\cal X}$}}
\newcommand{\bfcall}{\mbox{\boldmath ${\cal L}$}}
\newcommand{\bfcalf}{\mbox{\boldmath ${\cal F}$}}
\newcommand{\bfcalr}{\mbox{\boldmath ${\cal R}$}}
\newcommand{\bfcals}{\mbox{\boldmath ${\cal S}$}}
\newcommand{\bfcalw}{\mbox{\boldmath ${\cal W}$}}
\newcommand{\bfcalu}{\mbox{\boldmath ${\cal U}$}}
\newcommand{\bfcalv}{\mbox{\boldmath ${\cal V}$}}
\newcommand{\bfcalz}{\mbox{\boldmath ${\cal Z}$}}
\pagenumbering{roman}
\newcommand{\art}[6]{{\I[{\sc #1,}] {#2}, {\it #3}, {\bf #4}, {#5} {[#6]}}}
\newcommand{\ED}{\end{description}}
\newcommand{\I}{\item }
\newcommand{\ra}{\rm a}
\newcommand{\rb}{\rm b}
\newcommand{\rc}{\rm c}
\newcommand{\Hsp}{{\rm I}\!\!\,{\rm R}^n_+}
\newcommand{\Hsn}{{\rm I}\!\!\,{\rm R}^n_-}
\newcommand{\po}[1]{\mbox{$\displaystyle \frac{\mbox{$\partial #1$}}
{\mbox{$\partial x_{1}$}}$}}
\newcommand{\PO}[1]{\mbox{$\displaystyle \frac{\mbox{$\partial #1$}}
{\mbox{$\partial y_{1}$}}$}}
\newcommand{\OP}{\left(\Delta+2\lambda\PO{}\right)}
\newcommand{\op}{\left(\Delta+2\lambda\po{}\right)}
\newcommand{\ft}[1]{
\Frac{1}{(2\pi)^{n/2}}\Int{{\Bbb R}^{n}}{}e^{i{\bf x}\cdot \bfchi}
#1(\xi)d\xi}
\newcommand{\Ft}[1]{
\Frac{1}{2\pi}\Int{{\Bbb R}^{2}}{}e^{i{x}\cdot \xi}
#1(\xi)d\xi}
\newcommand{\Z}{\item[({\it a})]}
\newcommand{\B}{\item[({\it b})]}
\newcommand{\D}{\item[({\it d})]}
\newcommand{\E}{\item[({\it e})]}
\newcommand{\G}{\item[({\it g})]}
\def\tag{\renewcommand{\theequation}}
\newcommand{\Footnote}{~\footnote}
\newcommand{\ie}{{\it i.e.}}
\newcommand{\dist}{\mbox{\rm dist\,}}
\newcommand{\const}{\mbox{\rm const}}
\newcommand{\trace}{\mbox{\rm trace}}
\newcommand{\Bo}{\par\hfill{$\Box$}\par\noindent}
\newcommand{\Nor}[1]{\langle{#1}\rangle_q}
\newcommand{\vs}{\vspace*{.5cm}\par\noindent}
\newcommand{\Vs}{\vspace*{.6cm}\par\noindent}
\newcommand{\Vvs}{\vspace*{.7cm}\par\noindent}
\newcommand{\VVs}{\vspace*{.8cm}\par\noindent}
\newtheorem{definition}{Definition}[section]
\newcommand{\Bd}{\begin{definition}\begin{rm}}
\newcommand{\Ed}{\end{rm}\end{definition}}
\newtheorem{remark}{Remark}[section]
\newcommand{\Br}{\begin{remark}\begin{rm}}
\newcommand{\Er}{\end{rm}\end{remark}}
\newtheorem{proposition}{Proposition}[section]
\newcommand{\Bp}{\begin{proposition}\begin{sl}}
\newcommand{\EP}[1]{\end{sl}\label{proposition:#1}\end{proposition}}
\newcommand{\propref}[1]{{\rm Proposition \ref{proposition:#1}}}
\newcommand{\Bt}{\begin{theorem}\begin{sl}}
\newcommand{\Et}{\end{sl}\end{theorem}}
\newcommand{\Bl}{\begin{lemma}\begin{sl}}
\newcommand{\El}{\end{sl}\end{lemma}}
\newtheorem{theorem}{Theorem}[section]
\newtheorem{lemma}{Lemma}[section]
\newtheorem{lemmaA}{Lemma A.}
\newtheorem{corollary}{Corollary}[section]
\renewcommand{\eqref}[1]{{\rm (\ref{eq:#1})}}
\newcommand{\Bc}{\begin{corollary}\begin{sl}}
\newcommand{\Ec}{\end{sl}\end{corollary}}
\newcommand{\ET}[1]{\end{sl}\label{theorem:#1}\end{theorem}}
\newcommand{\EDD}[1]{\end{rm}\label{definition:#1}\end{definition}}
\newcommand{\EL}[1]{\end{sl}\label{lemma:#1}\end{lemma}}
\newcommand{\theoref}[1]{{\rm Theorem \ref{theorem:#1}}}
\newcommand{\ER}[1]{\end{rm}\label{remark:#1}\end{remark}}
\newcommand{\EC}[1]{\end{sl}\label{corollary:#1}\end{corollary}}
\newcommand{\defref}[1]{{\rm Definition \ref{definition:#1}}}
\newcommand{\remref}[1]{{\rm Remark \ref{remark:#1}}}
\newcommand{\cororef}[1]{{\rm Corollary \ref{corollary:#1}}}
\newcommand{\lemmref}[1]{{\rm Lemma \ref{lemma:#1}}}
\newcommand{\essup}[1]{{\rm ess}\,{{\displaystyle \sup_{\hspace*{-5mm}{#1}}}}}
\newcommand{\rosso}{\color{red}}
\newcommand{\nero}{\color{black}}


\pagenumbering{arabic}
\newcommand{\QED}{{\par\hfill$\square$\par}}
\renewcommand{\thefootnote}{(\arabic{footnote})}
\title{Global Weak Solutions to a Time-Periodic\\ Body-Liquid Interaction Problem} 
\author{Denis Bonheure~\thanks{D\'epartement de Math\'ematique,  Universit\'e Libre de Bruxelles, Belgium.
}\ \ \,\&\,\ Giovanni P. Galdi~\thanks{Department of Mechanical Engineering and Materials Science, University of Pittsburgh, USA. 
}
}
\date{\today}
\maketitle
\begin{abstract} We prove existence of time-periodic weak solutions to the coupled liquid-structure  problem constituted by an incompressible Navier-Stokes fluid interacting with a rigid body of finite size, subject to an {\em undamped} linear restoring force. The  fluid flow is generated by a uniform, time-periodic velocity field $\bfV$ far from the body. {We emphasize that our result is global, in the sense that no restriction is imposed on the magnitude of $\bfV$ and, rather remarkably, the frequency of $\bfV$ is entirely  arbitrary. Thus, in particular, it can coincide with any multiple of a natural frequency of vibration of the body so that, with this model, resonance cannot occur. Although based on the classical ``invading domains" technique, our approach requires several new ideas.} Indeed,  due to lack of sufficient dissipation, it appears quite unfeasible to show the existence of a fixed point of the Poincar\'e map at the finite-dimensional level along the Galerkin approximant. Therefore, unlike the usual strategy, such a result must be proven directly in a class of weak solutions, and therefore in the infinite-dimensional framework.\\[4mm] 
 \end{abstract}
\renewcommand{\theequation}{\arabic{section}.\arabic{equation}}
\tableofcontents  
\section*{Introduction} 
One among the many, captivating problems concerning the interaction of a liquid with an elastic structure arises when the liquid is in a time-periodic regime so as to generate a similar motion of the structure. Of particular interest is the case when the frequency of the flow, $\omega$, approaches or even coincides with a multiple of a natural frequency of vibration, $\omega_{\sf n}$, of the body. Then, the interaction may result in a resonance phenomenon, which is considered as primary responsible for possible failure of the structure \cite{Bev}. Such a problem falls in the general area of vibration-induced oscillations that has all along constituted a main focus of applied science and is at the heart of a vast engineering literature; see the monographs \cite{Bev, Dyr, Pai} and the bibliography therein.

In the classical  model employed  to investigate this type of questions, one regards the structure, $\mathscr B$,  as a rigid body subject to a linear restoring force, while the liquid, $\mathscr L$, is described by the Navier-Stokes equations  \cite{Bla,CM,MC,Will}. In such a framework, the occurrence of resonance is usually explained as follows \cite[Section 3.5]{Bev}, \cite{BH}. In the time-periodic regime of the liquid, the motion of the structure is that of a simple harmonic oscillator subject to a forced time-periodic motion caused by the action of the force, ${\sf F}$, exerted by $\mathscr L$ on $\mathscr B$. Then,  one can  show that for $\omega$   close to  (a multiple of) $\omega_{\sf n}$,  the amplitude of the forced oscillations of $\mathscr B$ becomes very large and, in absence of a damping mechanism, tends to infinity. A significant consequence of this argument is that, {\em without damping},  the coupled system $\mathscr S:=\mathscr B$--$\mathscr L$  cannot perform a time-periodic motion of {\em arbitrary} frequency $\omega$. While seemingly reasonable, this way of thinking appears rather simplistic. First, it decouples the highly coupled system $\mathscr S$, by somehow prescribing the action of $\mathscr L$ on $\mathscr B$. 
Moreover, it disregards the fact that  ${\sf F}$ splits into two competing components: the one that, indeed, forces the oscillations of $\mathscr B$, and the other, due to viscosity effects, which tends to absorb them. As a result, if the latter prevails over the former, occurrence of resonance can be excluded, even in {\em absence} of a structural damping mechanism. 
\par
Also motivated by the above considerations, very recently the authors and their collaborators started a {\em rigorous} and systematic  analysis of this model for different flow geometries, by investigating general mathematical properties, such as existence, uniqueness, stability and bifurcation of solutions   \cite{BBGGP,BG,BGG,BoGaGa,BoHiPaSpe,GazP,GPP,Patri}.   For other contributions related to a similar problematic, we also refer to \cite{GMZZ,Seb}. 

In this paper we aim at furnishing a further contribution, which may help to a better understanding of the resonance phenomenon. More specifically, we consider the general case when $\mathscr B$ (of arbitrary shape) is subject to a --possibly anisotropic-- linear restoring force, and is immersed in a Navier-Stokes liquid filling the whole space, $\Omega$, outside $\mathscr B$. The motion of the coupled system $\mathscr S$ is driven by a time-periodic uniform velocity, $\bfV=\bfV(t)$, of period $T$ impressed on the liquid at large distance from $\mathscr B$. We assume the worse-case scenario, namely,  {\em no external damping} mechanism acts on the structure, so that the only dissipative effect is due to the viscosity of $\mathscr L$. We then ask the question: Will the coupled system $\mathscr S$ perform a time-periodic motion of period $T$, for {\em arbitrary} $T$ and {\em arbitrary} magnitude of $\bfV$? A positive answer would suggest that the model used so far is probably not appropriate for  resonance studies, and that other aspects should be accounted for.    

The main achievement of this paper is to show that, indeed, the above question admits an affirmative answer, provided only that $\bfV$ has a mild degree of regularity and $\mathscr B$ is of class $C^2$. As expected, since we want to keep the ``size" of $\bfV$ arbitrary, we are lead to perform this study in a suitable class of weak solutions; see \defref{5.1}. The approach we use is, in principle, rather usual for time-periodic flow in exterior domains, and employs the ``invading domains" technique \cite{GK,GaSi2}. However, in the case at hand, its implementation  is by no means straightforward, and presents a number of difficulties that are described next.
 
We recall that the above technique develops along the following steps. One picks an increasing sequence of bounded domains, $\{\Omega_n\}$,  whose union coincides with $\Omega$, and suitably reformulates the original problem in each $\Omega_n$. Then, following an idea due to {\sc Prouse} \cite{Prouse}, a time-periodic solution in $\Omega_n$ is searched via the finite-dimensional Galerkin approximation, by showing with the help of Brouwer's theorem that the Poincar\'e map, ${\sf M}$, bringing  initial conditions into corresponding solutions at time $T$, has a fixed point. This fact, in conjunction with appropriate uniform estimates in $n$, allows one to construct a solution on each $\Omega_n$, and  eventually pass to the limit $n\to\infty$ to obtain a (weak) solution to the original problem. It should be emphasized that, for this approach to work, it is crucial that along the Galerkin approximant the  total energy of the system, $E$, in absence of forcing term, be bounded by an  exponentially decreasing function of time, to guarantee that ${\sf M}$ is a self-map. This means that the energy equation should contain a dissipative term proportional to $E$. It is exactly here that, in our case, the first difficulty arises. Actually, $E$ is the sum of the kinetic energy ($K_{\mathscr L}$) of the liquid, and the kinetic ($K_{\mathscr B}$) and potential ($U_{\mathscr B}$) energies of the body. By using viscosity dissipation combined with trace theorems we  obtain some damping terms for both $K_{\mathscr L}$ and $K_{\mathscr B}$ but, unsurprisingly, no damping for $U_{\mathscr B}$, namely, for the  oscillations amplitude of $\mathscr B$. In order to get the latter, one may think of  adapting a procedure seemingly introduced by {\sc Haraux} in the context of nonlinear wave equations \cite[p. 162 and $f\!f$]{HX};  see also \cite{HX2}. However, though this procedure  formally works on the original system of equations in $\Omega_n$, it is rather doubtful that it can be applied at the finite-dimensional level along the Galerkin approximant.   
 Therefore, we are forced to find a fixed point of the Poincar\'e map ${\sf M}$ directly on that system. Since we want existence for data of arbitrary size, the fixed point should be found in a class of weak solutions. Moreover, being now in an infinite-dimensional framework, ${\sf M}$ should also possess suitable compactness properties. This type of questions has been addressed by {\sc Prodi}  for classical Navier-Stokes equations \cite{Prodi}. Nevertheless, his method requires uniqueness and continuous dependence of solutions upon the data in the energy norm,  properties that, to date, are only known to hold in two dimensions.  As a consequence,  {\sc Prodi}'s method is inapplicable to our case, which produces yet another  difficulty. To overcome this issue, we introduce an entirely different strategy that combines a suitable mollification of the nonlinear term in the original problem in $\Omega_n$, along with the use of time-weighted norms that vanish at $t=0$. Incidentally, we remark that, by using  this new strategy,  one could extend {\sc Prodi}'s result to arbitrary dimension $d\ge2$.

Our approach employs the following steps. Thanks to the regularization procedure,  we are able to show the existence of global weak solutions to the mollified problem for initial data possessing only finite energy, that however are strong (\`a la Prodi-Ladyzhenskaya) at any  positive time; see \lemmref{3.3}.  In this class of solutions it is easily shown that the Poincar\'e map is compact in the energy space. Thus, in order to prove existence of time-periodic solutions to the mollified problem, it remains to ascertain that ${\sf M}$ maps some ball in the energy space into itself. This property is  shown by adapting {\sc Haraux}'s argument mentioned earlier on, which gives the desired dissipation also for $U_{\mathscr B}$; see \lemmref{4.3}. In this way, by letting the mollification parameter go to $0$, we finally deduce the existence of a time-periodic weak solution in every $\Omega_n$; see \propref{4.1}.  The last step is to let $n\to\infty$,   and this brings the last difficulty.  Actually, for the dissipative term for $U_{\mathscr B}$, it  appears  not possible to provide an estimate that is {\em uniform} in $n$, which means that, as $n\to\infty$, we have no control on the amplitude of the oscillations of $\mathscr B$. However, we prove a uniform control on the velocity of the center of mass of $\mathscr B$ along with estimates for the time derivative of the flow velocity in suitable distributional spaces.  Combining these estimates   allows us to deduce, by a (local) compactness argument, existence of time-periodic weak solutions for the original problem in the whole $\Omega$; see Section \ref{sec:Fin}.
\par
We conclude this introductory section with the following remark. In the model considered in the present paper (and precisely presented in the next section), the body $\mathscr B$ can move only by translational motion and is not free to rotate. Removing this constraint and adding a corresponding restoring torque  would lead to a more complete model that might present resonance phenomena, owing to possible  interaction of the different degrees of freedom.  Plainly, the mathematical analysis of such a model becomes  far more complicated, also because of the notorious difficulty due to the presence of an unbounded coefficient in the flow equations \cite{Gah}.  The investigation of this topic   will be the object of future work. 
\par
The plan of the paper is as follows. After introducing the mathematical formulation of the problem in Section 1, in the following Section 2 we introduce the basic function spaces and collect some preliminary results that, among other things, help us furnishing a suitable reformulation of the original problem. In Subsection 3.1 we state our main result in \theoref{5.1} and present the strategy  we use for its proof. Subsection 3.2--3.5 are dedicated to the proof of existence of time-periodic solutions to a suitable modification of the problem in an arbitrary bounded domain, strictly containing $\mathscr B$, along with uniform estimates independent of the ``size" of the domain. With this result in hand, in the final Subsection 3.6 we are able to produce a full proof of \theoref{5.1}.

\setcounter{equation}{0} 
\section{Formulation of the Problem}\label{sec:Form}
Consider a rigid body $\mathscr B$, occupying the closure of the  bounded domain $\Omega_0$, completely surrounded by a Navier-Stokes liquid, $\mathscr L$,  filling  the entire space, $\Omega$, outside $\mathscr B$. $\mathscr B$ is subject to an elastic restoring force, $\sf{\bf R}$,  applied at its center of mass $G$. We take $\sf{\bf R}$ to be linear, but not necessarily isotropic, that is,
$$
\sf{\bf R}=-\hat{\mathbb A}\cdot\bfdelta\,,
$$
where $\bfdelta:=\vec{GO}$,  $O$ a fixed point, and $\hat{\mathbb A}$ a $3\times3$ symmetric and positive definite real matrix (stiffness matrix).
We assume that $\mathscr B$ can move only by translational motion, which can be accomplished by having a suitable torque acting on it.  
The motion of the coupled system body-liquid is driven by a time-periodic flow of $\mathscr L$ imposed at ``large" spatial distances from $\mathscr B$ and characterized by a uniform $T$-periodic velocity field $-\hat{\bfV}$, where $\hat{\bfV}$ is a bounded function of time $t$ only, satisfying therefore
 $$\hat{\bfV}(t)=\hat{\bfV}(t+T)\,,\ \ \mbox{ for all $t\in\real$}\,.
$$    
Denote by $L$ and $M$ diameter and mass of $\mathscr B$, and by $\rho$ and $\mu$ density and shear viscosity coefficient of the liquid. Then,
the $T$-periodic motion of the coupled system body-liquid when referred to a body-fixed frame $\calf\equiv\{G,\bfe_i\}$  is governed by the following set of dimensionless equations, see for instance \cite[Section 1]{Gah},
   
\be\left\{\ba{l}\medskip\left.\ba{r}\medskip
\partial_t\bfv+\lambda(\bfv-{\bfgamma})\cdot\nabla\bfv=\Delta\bfv-\nabla {\sf p}\\
\Div\bfv=0\ea\right\}\ \ \mbox{in $\Omega\times\real$}\,,\\ \medskip
\ \ \bfv(x,t)={\bfgamma}(t)\,, \ \mbox{ $(x,t)\in\partial\Omega\times\real$}\,;\ \\ \medskip
\ \ \Lim{|x|\to\infty}\bfv(x,t)=-{\bfV}(t)\,,\ t\in\real\,,\\
\medskip\left.\ba{r}\medskip
\dot{\bfgamma}+\mathbb A\cdot\bfdelta+\varpi\Int{\partial\Omega}{} \mathbb T(\bfv,{\sf p})
\cdot\bfn=\0 \\
\dot{\bfdelta}=\bfgamma
\ea\right\}\ \ \mbox{in $\real$.}
\ea\right.
\eeq{01}
Here, $\bfv$ and ${\sf p}$ are (non-dimensional) velocity and pressure fields of the liquid,   while
$$
 \mathbb A:=\frac{\rho^2 L^4}{M\mu^2}\hat{\mathbb A} \,,\ \ \varpi:=\frac{\rho L^3}{M}\,,\ \ \lambda:=\frac{\rho V_\infty L}{\mu}\,,\ \ \bfV:=\frac{\hat{\bfV}}{V_\infty}\,,\ \
V_\infty:=\sup_{t\in\real}|\hat{\bfV}(t)|\,.
$$
 
Moreover, 
$$
\mathbb T(\bfz,\psi):=2\,\mathbb D(\bfz)-\psi\,\mathbb I\,,\ \ \ \mathbb D(\bfz):=\half\big(\nabla\bfz+(\nabla\bfz)^\top\big)\,,
$$
with $\mathbb I$ identity matrix, is the (dimensionless) Cauchy stress tensor, and  $\bfn$ the unit outer normal at $\partial\Omega$.
\par
Notice that, with the above  non-dimensionalization, we have 
\be
\sup_{t\in \real}|\bfV(t)|=1.
\eeq{SV}
Our ultimate goal is to show that, for any given $\lambda, \mathbb A, \varpi$ and $T>0$, and any  (sufficiently smooth) $\bfV$, problem \eqref{01} has at least one,  suitably defined,  $T$-periodic weak solution $(\bfv,\bfgamma,\bfdelta)$.\par     
\setcounter{equation}{0} 
\section{Preliminaries}
\subsection{Functional Spaces and Some Related Properties}\label{notation}
Before describing our functional framework, we begin with some  notation. We indicate by $\Omega\subset \real^3$   the exterior domain of class $C^2$, defined as  complement of the closure of the  bounded domain $\Omega_0$  occupied by $\mathscr B$. 
We take the origin of coordinates in the interior of $\Omega_0$, and set  $B_r:=\{x\in\real^3:\,|x|<r\,,\ r>0\}$, and 
$$\Omega_R:=\Omega\cap B_R,\  R>R_*:={\rm diam}\,\Omega_0.$$
As customary, for $A$ a domain of $\real^3$,  $L^q=L^q(A)$, $W^{m,2}=W^{m,2}(A)$, $q\in[1,\infty]$, $m\in\nat$, are Lebesgue and  Sobolev spaces with norm  $\|\cdot\|_{q,A}$, and $\|\cdot\|_{m,2,A}$. By $(\,\ ,\,\ )_A$ we indicate the $L^2(A)$-scalar product.  Furthermore, $D^{m,q}=D^{m,q}(A)$ is the homogeneous Sobolev space with semi-norm $\sum_{|l|=m}\|D^lu\|_{q,A}$. 
%
In all the above notation we shall typically omit the subscript ``$A$", unless confusion arises. For a Banach space $X$, we may, occasionally, indicate its norm by $\|\cdot\|_X$. Finally, by $L^q(I;X)$, $W^{1,q}(I;X)$ and $C^m(I;X)$, where $I$ is a real interval, we shall denote the classical Bochner spaces.
\par
If ${\sf A}\subseteq \real^3$ is a domain with ${\sf A}\supset \bar{\Omega_0}$, let
$$\ba{ll}\medskip
\calk=\mathcal K({\sf A}):=\big\{\bfphi\in C_0^\infty({\sf A}) : 
\exists\,\hat{\bfphi}\in\real^3 \mbox{ s.t. }\bfphi(x)=\hat{\bfphi} \mbox{ in a neighborhood of }\Omega_0\big\}\,,
\\ \medskip
\calc=\mathcal C({\sf A}):=\{\bfphi\in\calk({\sf A}):\ \Div\bfphi=0\ \mbox{in ${\sf A}$}\}\,,\\
\calc_0=\mathcal C_0({\sf A}):=\{\bfphi\in\calc({\sf A}): \hat{\bfphi}=\0\}\,.
\ea
$$
In $\calk({\sf A})$ we introduce the scalar product
\be
\langle \bfphi,\bfpsi\rangle_{\sf A}:=
\varpi^{-1}\,\hat{\bfphi}\cdot\hat{\bfpsi}+(\bfphi,\bfpsi)_{{\sf A}\cap\Omega}\,,\ \ \bfphi,\bfpsi\in\calk\,,
\eeq{0.0}
and define 
\be\ba{ll}\medskip
\call^2(\real^3):= \,\big\{\mbox{completion of $\calk(\real^3)$ in the norm induced by \eqref{0.0}}\big\}\,,\\ \medskip
\calh(\real^3):=\,\big\{\mbox{completion of $\calc(\real^3)$ in the norm induced by \eqref{0.0}}\big\}\,,\\ 
\calg(\real^3):=\big\{\bfh\in \call^2(\real^3): \, \exists\, p\in D^{1,2}(\Omega) \mbox{ s.t. } \bfh=\nabla p \mbox{ in } \Omega,
\mbox{and}\ \bfh=-\varpi\int_{\partial\Omega}p\,\bfn   \ \mbox{in}\ \Omega_0\,\big\}\,.

\ea
\eeq{spazi}
It is shown in \cite[Theorem 3.1 and Lemma 3.2]{ALS} that
$$\ba{ll}\medskip
\call^2(\real^3)=\{\bfu\in L^2(\real^3): \ \bfu=\hat{\bfu}\ \mbox{in}\ \Omega_0, \ \mbox{for some $\hat{\bfu}\in \real^3$}\}\\
\calh(\real^3)=\{\bfu\in \call^2(\real^3): \ \Div\bfu=0\,\}\,,
\ea
$$
along with the following orthogonal decomposition \cite[Theorem 3.2]{ALS}
\be
\call^2(\real^3)=\calh(\real^3)\oplus\calg(\real^3)\,.
\eeq{Helm}
We next define the space
$$
\cald^{1,2}=\cald^{1,2}(\real^3):=\,\big\{\mbox{completion of $\calc(\real^3)$ in the norm $\|\mathbb D(\cdot)\|_2$}\big\}
$$
whose basic properties are collected in the next lemma; see \cite[Lemmas 9--11]{Gah}.
\Bl  ${\cald^{1,2}}$ is a separable Hilbert space when equipped with the scalar product
$$
(\mathbb D(\bfu_1),\mathbb D(\bfu_2))\,,\ \ \bfu_i\in {\cald^{1,2}}\,, \, \ i=1,2\,.
$$
Moreover, we have the characterization:
\be
{\cald^{1,2}}=\big\{\bfu\in L^6(\real^3)\cap D^{1,2}(\real^3)\,;\ \Div\bfu=0\,;\,
\bfu=\hat{\bfu} \ \mbox{in $\Omega_0$}\,,\ \mbox{for some $\hat{\bfu}\in\real^3$} \big\}\,.
\eeq{1.7777}
Also, for each $\bfu\in{\cald^{1,2}}$, it holds
\be
\|\nabla\bfu\|_2=\sqrt{2}\|\mathbb D(\bfu)\|_2\,,
\eeq{1.8}
and
\be
\|\bfu\|_6\le \kappa_0\,\|\mathbb D(\bfu)\|_2\,,
\eeq{1.9}
for some numerical constant $\kappa_0>0$. Finally, 
there is another positive constant $\kappa_1$ such that
\be
|\hat{\bfu}|\le \kappa_1\,\|\mathbb D(\bfu)\|_2\,.
\eeq{1.10}
\EL{1.1_1}
\par
Along with the spaces $\call^2,\calh$, and  $\cald^{1,2}$ defined above, we introduce suitable ``local" versions of these spaces. Precisely, we set
$$\ba{ll}\medskip
\call^2(B_R):= \{\bfphi\in L^2(B_R): \, \bfphi|_{\Omega_0}=\hat{\bfphi}\,\ \mbox{for some $\hat{\bfphi}\in\real^3$}\}\,,\\ \medskip
\calh(B_R):=\{\bfphi\in \call^2(B_R):\, \Div\bfphi=0\,,\ \bfphi\cdot\bfn|_{\partial B_R}=0\}\,,\\ \medskip
\cald^{1,2}(B_R):= \{\bfphi\in W^{1,2}(B_R): \, \Div\bfphi=0\,,\ \bfphi|_{\Omega_0}=\hat{\bfphi}\,\ \mbox{for some $\hat{\bfphi}\in\real^3$}\,,\ \bfphi|_{\partial B_R}=\0 \}\,,\\
\cald^{1,2}_0(\Omega_R):= \{\bfphi\in \cald^{1,2}(B_R): \, \hat{\bfphi}=\0 \}\,.
\ea
$$
Then $\calh(B_R)$ and $\cald^{1,2}(B_R)$ are Hilbert spaces with scalar products
$$
\langle\bfphi_1,\bfphi_2\rangle_{B_R}\,,\ \ \bfphi_i\in\calh(B_R)\,;\ \ \ (\mathbb D(\bfpsi_1),\mathbb D(\bfpsi_2))_{B_R}\,,\ \bfpsi_i\in \cald^{1,2}(B_R)\,,\ \ i=1,2.
$$
Moreover, the following decomposition holds, analogous to \eqref{Helm}  \cite[Theorem 3.1 and Lemma 3.2]{ALS}
\be
\call^2(B_R)=\calh(B_R)\oplus\calg(B_R)\,,
\eeq{Helm1}
where $\calg(B_R)$ is defined as in \eqref{spazi}$_3$, by replacing $\Omega$ with $\Omega_R$.
\par
Finally, the dual spaces of $\cald^{1,2}(B_R)$ and $\cald^{1,2}_0(\Omega_R)$ will be denoted by $\cald^{-1,2}(B_R)$ and $\cald^{-1,2}_0(\Omega_R)$, respectively.
\Br The space $\cald^{1,2}(B_R)$ can be viewed as a subspace of $W^{1,2}(\Omega)\cap\cald^{1,2}(\real^3)$, by extending its generic element to 0 in $\real^3\backslash B_R$. Therefore, all the properties mentioned in \lemmref{1.1_1} continue to hold for $\cald^{1,2}(B_R)$. 
\ER{2.1}

We next recall basic facts about the mollification of fields in $\cald^{1,2}(B_R)$.
Let $\bfu\in\cald^{1,2}(B_R)$ and continue to denote by $\bfu$ its extension to $W^{1,2}(\real^3)\cap\cald^{1,2}(\real^3)$, in the sense of \remref{2.1}. Moreover, let
$\eta_0>0$ be small enough so that the domain
\be
\Omega_{00}:=\{x\in\Omega_0:\ \dist(x,\partial\Omega)>\eta_0\} 
\eeq{Om0}
is not empty.
For a given $\eta\in (0,\eta_0)$,  we then indicate by $\bfu_\eta$ the (Friederichs) mollifier of $\bfu$, namely,  
\be
\bfu_\eta(x):=\Int{\real^3}{}k_\eta(x-y)\bfu(y){\rm d}y\,,
\eeq{regu}
where 
$$
k_\eta(\xi):=\eta^{-3}k(\xi/\eta)\,;\ \  k\in C_0^\infty(B_1)\,,\ \ \int_{\real^3}k(x){\rm d}x=1\,.
$$
\Bl Let $\bfu\in\cald^{1,2}(B_R)$ and $\bfu_\eta$ be defined by \eqref{regu} with $\eta\in (0,\eta_0)$. The following properties hold:
\begin{itemize}
  \item [{\rm (a)}] $\Div\bfu_\eta(x)=0$, for all $x\in\real^3$\,;
\item [{\rm (b)}] $\bfu_\eta(x)=\hat{\bfu}$, for all $x\in \bar{\Omega_{00}}$\,,\ \mbox{ and with $\hat{\bfu}=\bfu|_{\Omega_0}$}\,;  
\item [{\rm (c)}] for any $\bfw\in\cald^{1,2}(B_R)$, 
 $$
\int_{\Omega_R}(\bfu_\eta-\hat{\bfu})\cdot\nabla\bfw\cdot\bfw=0\,.
$$
\end{itemize}
\EL{3.2}
{\em Proof.} Since $\Div\bfu(x)=0$, $x\in\real^3$ and $\bfu(x)=\hat{\bfu}$, $x\in \Omega_0$, both (a) and (b) follow from the properties of mollifiers and the definition of $\Omega_{00}$ in \eqref{Om0}.  Integrating by parts and using $\bfw|_{\partial B_R}=\0$, $\bfw|_{\partial \Omega}=\hat{\bfw}$ along with (a), we get
$$
\int_{\Omega_R}(\bfu_\eta-\hat{\bfu})\cdot\nabla\bfw\cdot\bfw=\half|\hat{\bfw}|^2\int_{\partial\Omega}(\bfu_\eta-\hat{\bfu})\cdot\bfn=\half|\hat{\bfw}|^2\int_{\partial\Omega}\bfu_\eta\cdot\bfn=\half|\hat{\bfw}|^2\int_{\Omega_0}\Div\bfu_\eta=0\,.
$$
\hfill$\square$\par
\subsection{Reformulation of the Problem with Velocity Fields Vanishing at Infinity}\label{reform}
In order to solve Problem \eqref{01} formulated in Section \ref{sec:Form}, it is convenient and customary to deal with velocity fields vanishing  as $|x|\to\infty$.   This requires an appropriate lifting of $\bfV$, which will be accomplished with the help of the following result.
\Bl Let $\eta_0>0$, $\Omega_{00}$ be as in \eqref{Om0},
and set
$$
\Omega^0:=\real^3\backslash\bar{\Omega_{00}}\,.
$$
Then, for any $\varepsilon>2/\ln(1/\eta_0)$, there exists $\bfU=\bfU(\varepsilon;x,t)$, $(x,t)\in \Omega^0\times \real$ such that
\begin{itemize}
  \item [{\rm (i)}] $\bfU(t)\in C^\infty(\Omega^0)\,,\ \ t\in\real$\,;  
  \item [{\rm (ii)}] $\Div\bfU(x,t)=0\ \ (x,t)\in\Omega^0\times\real$\,;
   \item [{\rm (iii)}]  $\supp(\bfU(t))\subset \{x\in\Omega^0:\, 0\le \dist(x,\partial\Omega_{00})\le {\rm e}^{-1/\varepsilon}\}=:\Omega_\varepsilon\,,$ \ for all $t\in \real$;
 \item [{\rm (iv)}] $\bfU(x,t)=\bfV(t)\,,\ \ (x,t)\in \Omega_{\frac{\varepsilon}2}\times\real\supset \partial\Omega\times\real\,;$  
 \item [{\rm (v)}]  $\|\bfU(t)\|_{2,2}\le c_\varepsilon\,|\bfV(t)|$\,;
   
  \item [{\rm (vi)}]  $\|\partial_t\bfU(t)\|_{2}\le c_\varepsilon\,|\dot \bfV(t)|$\,;
   
\item [{\rm (vii)}] If $\hat{\Omega}\supset \Omega_\varepsilon$, and $\bfw,\bfz\in W^{1,2}(\hat{\Omega})$ with $\bfw|_{\partial\Omega_{00}}=\0$, then
$$
\int_{\hat{\Omega}}\left|\bfw\cdot\nabla\bfz\cdot\bfU\right|\le c_0\,\varepsilon \|\nabla\bfw\|_2\|\nabla\bfz\|_2\,,
$$
where $c_0$ is a positive constant independent of $\varepsilon$.
\end{itemize}
\EL{1.1}
{\em Proof.} Let $\psi=\psi(r)$, $r\in [0,\infty)$ be a smooth, non-decreasing real function such that $\psi(r)=0$ if $r\le1$, and $\psi(r)=1$, if $r\ge2$, and set
$$
\phi(\varepsilon;x)=\psi(-\varepsilon \ln d(x)),
$$
where $d(x):=\dist(x,\partial\Omega_{00})$.
Clearly,
$$
\phi(\varepsilon;x)=\left\{\ba{ll}\medskip 1 & \mbox{if}\ d(x)\le {\rm e}^{-2/\varepsilon}\\
0 & \mbox{if}\ d(x)\ge {\rm e}^{-1/\varepsilon}
\ea\right.\,,
$$
and, moreover,
\be
|\nabla\phi(\varepsilon;x)|\le \frac{c\,\varepsilon}{d(x)}\,,
\eeq{1.2}
with $c$ independent of $\varepsilon$. Let
$$
\bsfU(x,t):=x_3V_2(t)\bfe_1+x_1V_3(t)\bfe_2+x_2V_1(t)\bfe_3
$$
and define
$$
\bfU(\varepsilon;x,t):=\curl(\phi(\varepsilon;x)\bsfU(x,t))\,.
$$
Since $\curl\bsfU=\bfV$, we get
\be
\bfU(\varepsilon;x,t)=\phi(\varepsilon;x)\bfV(t)-\bsfU(x,t)\times\nabla\phi(\varepsilon;x)\,,
\eeq{1.3}
which, by the properties of $\phi$ and the choice of $\varepsilon$, shows the validity of (i)--(vi) above. Moreover, by Schwarz inequality,  \eqref{1.3}, \eqref{1.2}, and observing that $\sup_{\Omega_\varepsilon}|\bsfU|\le c\,|\bfV(t)|$,
we deduce~\footnote{Recall \eqref{SV}.\label{foot:1}}
\be\ba{rl}\medskip
\left(\Int{\hat{\Omega}}{}\left|\bfw\cdot\nabla\bfz\cdot\bfU\right|\right)^2= \left(\Int{{\Omega_\varepsilon}}{}\left|\bfw\cdot\nabla\bfz\cdot\bfU\right|\right)^2 & \le \|\nabla\bfz\|_2^2\Int{{\Omega_\varepsilon}}{}|\bfU|^2|\bfw|^2\\ \medskip
& \le c\,\|\nabla\bfz\|_2^2\Int{{\Omega_\varepsilon}}{}\left(\phi^2|\bfw|^2+\varepsilon^2\Frac{|\bfw|^2}{d^2}\right)\\ \medskip
&\le c\|\nabla\bfz\|_2^2\Int{{\Omega_\varepsilon}}{}\left(d^2+\varepsilon^2\right)\Frac{|\bfw|^2}{d^2}\\ & \le c\,\varepsilon^2\,\|\nabla\bfz\|_2^2\Int{{\Omega_\varepsilon}}{}\Frac{|\bfw|^2}{d^2}\,,
\ea
\eeq{1.4} 
where $c$ is independent of $\varepsilon$.
Using Hardy inequality, see e.g. \cite[Lemma III.6.3]{Gab}, we conclude that 
$$
\left(\Int{\hat{\Omega}}{}\left|\bfw\cdot\nabla\bfz\cdot\bfU\right|\right)^2 
\le c\,\varepsilon^2\,\|\nabla\bfz\|_2^2\Int{{\Omega_\varepsilon}}{}\Frac{|\bfw|^2}{d^2} 
\le c\,\varepsilon^2\,\|\nabla\bfz\|_2^2\Int{{\Omega_\varepsilon}}{}|\nabla\bfw|^2,
$$
and the property (vii) follows. 
\par\hfill$\square$\par
With \lemmref{1.1} at hand, we now rewrite \eqref{01} with the new unknown velocity field
$$
\bfu=\bfv+\bfV-\bfU
$$
so that $(\bfu,\bfgamma,\bfdelta)$ now satisfies the following problem with source terms:

\be\left\{\ba{l}\medskip\left.\ba{r}\medskip
\partial_t\bfu+\lambda\left[(\bfu-{\bfgamma}+\bfU-\bfV)\cdot\nabla\bfu+(\bfu-\bfgamma)\cdot\nabla\bfU\right]=\Div\mathbb T(\bfu,{p})+\bff\\
\Div\bfu=0\ea\right\}\ \ \mbox{in $\Omega\times\real$}\,,\\ \medskip
\ \ \bfu(x,t)={\bfgamma}(t)\,, \ \mbox{ $(x,t)\in\partial\Omega\times\real$}\,;\\ \medskip
\ \ \Lim{|x|\to\infty}\bfu(x,t)=\0\,,\ t\in\real\,,\\
\medskip\left.\ba{r}\medskip
\dot{\bfgamma}+\mathbb A\cdot\bfdelta+\varpi\Int{\partial\Omega}{} \mathbb T(\bfu,p)\cdot\bfn=\bfF\\
\dot{\bfdelta}=\bfgamma
\ea\right\}\ \ \mbox{in $\real$.}
\ea\right.
\eeq{1.5}
where
\be\ba{ll}\medskip
\bff:=(\bfV-\bfU)\cdot\nabla\bfU+\Delta\bfU-\partial_t\bfU\,,\\ \medskip
\bfF:=-\varpi\,{\rm vol}\,(\mathscr B)\,\dot{\bfV}\,,\\ \medskip p:={\sf p}-\dot{\bfV}\cdot\bfx\,.
\ea
\eeq{1.6}
Observe that, by the properties of $\bfU$ listed in \lemmref{1.1}, and \eqref{SV},  we have, in particular, that
\be
\supp(\bff(t))\subset\Omega_\varepsilon\,, \ \ \mbox{for all $t\in\real$}\,;\ \ \|\bff(t)\|_2\le c\,(|\bfV(t)|+|\dot{\bfV}(t)|)\,. 
\eeq{1.7}
We thus obtain that the original problem \eqref{01} has been formally and equivalently reformulated as \eqref{1.5}--\eqref{1.6} where the velocity of $\mathscr L$ vanishes at infinity,  whereas  non-zero prescribed external $T$-periodic forces are now acting on both $\mathscr L$ and $\mathscr B$. 
\par

\setcounter{equation}{0} 
\section{Existence of $T$-periodic Weak Solutions}
\subsection{Definition of $T$-periodic Weak Solutions}

We first need to introduce a suitable general class of $T$-periodic test functions. Precisely,
let $A$ be either $B_R$ or $\real^3$. By
$\calc_\sharp(A)$, we denote  the space of restriction to $[0,T]$ of functions $\bfphi \in C^{1}(A\times \mathbb{R})$, satisfying: 
\begin{itemize}
\item[{\rm (a)}] $\Div  \bfphi(x,t) = 0$ for $(x,t)\in A \times \mathbb{R}$\,; 
\item[{\rm (b)}] $\bfphi (x,t) = \hat\bfphi (t)$, some $\hat\bfphi \in C^{1}(\mathbb{R})$,   for $x$ in a neighborhood 
of $\Omega_0$ and $t \in \mathbb{R}$\,;
 
\item[{\rm (c)}] $\supp_x{\bfphi(x,t)}\subset A$ for all $t \in\real$\,;
\item[{\rm (d)}] $\bfphi(x,t+T)=\bfphi(x,t)$ for all $(x,t) \in A\times\mathbb{R}$\,.
\end{itemize}

We are then able to give the definition of $T$-periodic weak solution to \eqref{1.5}. Testing \eqref{1.5}$_1$ by arbitrary $\bfphi\in \calc_\sharp(\real^3)$, integrating by parts over $\Omega\times [0,T]$ and employing \eqref{1.5}$_{2-5}$, we infer 
\[
\begin{split}
\langle\bfu(T),\bfphi(T)\rangle-\langle\bfu(0),\bfphi(0)\rangle
= &
-\Int0T\big[-\langle\bfu,\partial_t\bfphi\rangle+\lambda\left((\bfu-{\bfgamma}+\bfU-\bfV)\cdot\nabla\bfu+(\bfu-\bfgamma)\cdot\nabla\bfU,\bfphi\right)\\ 
& + \left.2(\mathbb D(\bfu),\mathbb D ({\bfphi}))+\varpi^{-1}  \hat\bfphi \cdot\mathbb A\cdot\bfdelta-(\bff,\bfphi)-\bfF\cdot  \hat\bfphi \right]{\rm d}t\,, 
\end{split}
\]
where we recall that $(\cdot,\cdot)\equiv (\cdot,\cdot)_{\Omega}$ is the $L^2(\Omega)$-scalar product whereas 
$$
\langle \bfphi,\bfpsi\rangle=
\varpi^{-1}\,\hat{\bfphi}\cdot\hat{\bfpsi}+(\bfphi,\bfpsi)_{\Omega}.$$
Furthermore, from \eqref{1.5}$_6$, we also have
$$
\bfdelta(T)-\bfdelta(0)=\int_0^T\bfgamma(t)\,{\rm d}t\,.
$$
Thus, if $(\bfu,\bfgamma,\bfdelta)$ is a $T$-periodic (sufficiently smooth) solution to \eqref{3.1}, 
then $\bfgamma$ has zero average. Since the periodicity also implies
$$\langle\bfu(T),\bfphi(T)\rangle=\langle\bfu(0),\bfphi(0)\rangle,$$
we deduce that any such $T$-periodic solution satisfies
\be
\begin{split}
\Int0T\Big[-\langle\bfu,\partial_t\bfphi\rangle+\lambda\big((\bfu-{\bfgamma}+\bfU-\bfV)\cdot\nabla\bfu & + (\bfu-\bfgamma)\cdot\nabla\bfU,\bfphi)+2(\mathbb D(\bfu),\mathbb D ({\bfphi})\big)\\  & +\varpi^{-1}\hat\bfphi\cdot\mathbb A\cdot\bfdelta -(\bff,\bfphi)-\bfF\cdot \hat\bfphi\Big]{\rm d}t=0\,,\\ 
& \dot{\bfdelta}=\bfgamma\,,\ \ \Int0T\bfgamma(t){\rm d}t=\0,
\end{split}
\eeq{4.521}
for arbitrary test function $\bfphi$  in $\calc_\sharp(\real^3)$.
Conversely,  with the help of the decomposition \eqref{Helm}, it is easy to see
that every sufficiently smooth functions $(\bfu,\bfgamma,\bfdelta)$ obeying \eqref{4.521}, is a $T$-periodic solution to \eqref{1.5}. With this in mind, we give the following definition of weak solution.
\Bd The triple $(\bfu,\bfgamma,\bfdelta)$ is a $T$-\emph{periodic weak solution} to \eqref{1.5} if
\begin{itemize} 
\item [{\rm (i)}] $\bfu \in L^2(0,T;{\mathcal D}^{1,2}(\real^3))$\,, with $\bfu(x,t)|_{\partial\Omega}=\bfgamma(t)$\,, a.a. $t\in [0,T]$\,, 
$\bfgamma\in L^2(0,T;\mathbb{R}^3)$\,; 
\item [{\rm (ii)}] $\bfdelta \in W^{1,2}(0,T;\real^3)$\,;
\item [{\rm (iii)}] $(\bfu,\bfgamma,\bfdelta)$ satisfies \eqref{4.521}.
\end{itemize}
\EDD{5.1}
\Br In view of   \lemmref{1.1_1} , it is easy to check that the  integral in \eqref{4.521} is well defined for a weak solution. Likewise, the boundary condition in {\rm (i)} is meaningful in the trace sense.
\ER{4.3}
\par

\subsection{Statement of the Main Theorem and Strategy of the Proof}

The main contribution of this paper is expressed by the following result.
\Bt
Suppose $\bfV\in W^{1.2}(0,T;\real^3)$, be $T$-periodic for some $T>0$. Then,  there exists  at least one corresponding $T$-periodic weak solution to \eqref{1.5}. This solution satisfies the estimate
\be
\int_0^T(\|\nabla\bfu(t)\|_{2}^2+|\bfgamma(t)|^2){\rm d}t\le C\,\int_0^T(|\bfV|^2+|\dot{\bfV}|^2){\rm d}t\,,
\eeq{5.2}
where the constant $C$ depends only on $\Omega$ and the physical parameters of the body and the liquid.\label{T5.1}
\ET{5.1}

Our strategy to prove \theoref{5.1}  goes as follows. In a first step, we shall consider a suitable modification of problem  \eqref{1.5}--\eqref{1.6} in a generic bounded domain of the type $\Omega_R$.  Precisely, for any $R>3R_*$  and such that $B_R\supset \Omega_\varepsilon$,\footnote{Later on, the parameter $\varepsilon$ will be fixed in terms of the data, so that this request is meaningful; see \eqref{eta}.}  we shall prove that the problem
\be\left\{\ba{l}\medskip\left.\ba{r}\medskip
\partial_t\bfu+\lambda\left[(\bfu-{\bfgamma}+\bfU-\bfV)\cdot\nabla\bfu+(\bfu-\bfgamma)\cdot\nabla\bfU\right]=\Div\mathbb T(\bfu,{p})+\bff\\
 \Div\bfu=0\ea\right\}\ \ \mbox{in $\Omega_R\times\real$}\,,\\ \medskip
\ \ \bfu(x,t)={\bfgamma}(t)\,, \ \mbox{ $(x,t)\in\partial\Omega\times\real$}\,,\\ \medskip
\ \ \bfu(x,t)=\0\,, \ \mbox{ $(x,t)\in\partial B_R\times\real$},\\
\left.\ba{r}\medskip
\dot{\bfgamma}+\mathbb A\cdot\bfdelta+\varpi\Int{\partial\Omega}{} \mathbb T(\bfu,p)\cdot\bfn=\bfF\\
\dot{\bfdelta}=\bfgamma
\ea\right\}\ \ \mbox{in $\real$}\,,
\ea\right.
\eeq{3.1}
with $\bff$ and $\bfF$  given in \eqref{1.6} has at least one $T$-periodic (weak) solution $(\bfu,\bfgamma,\bfdelta)$ that, in addition, obeys certain bounds in terms of the data,  uniformly with respect  to $R$. We then let $R\to \infty$ along a sequence and prove that the corresponding solutions will converge to a weak solution to the original problem. 
\par

Even though this approach is classical, its implementation in the present setting is by no means straightforward, due to the fact that the spring has no damping. In particular, the classical method of showing the existence of a fixed point for the Poincar\'e map at the finite-dimensional level (along the Galerkin approximations) fails, due to the lack of ``sufficient dissipation." We are thus lead to prove this existence at the infinite-dimensional level.  

We develop the steps towards the proof of \theoref{5.1} in the next subsections. 
We first show existence and uniqueness of strong solutions to the initial-value problem associated to  a regularized version of \eqref{3.1}, where the nonlinear term has been suitably mollified. Then,  we prove that the Poincar\'e map associated to this problem has a fixed point, which leads to the existence of $T$-periodic strong solutions. Finally, we let the mollifying parameter to $0$, thus obtaining the same result for the original problem \eqref{3.1} in the class of weak solutions. 

The last step is to let $R\to \infty$.  We are not able to provide uniform estimates on the dissipation in $\bfdelta$, which means that, as $R\to\infty$, we have no control on the amplitude of the oscillations. To overcome this issue, we prove the velocity $\dot\bfdelta$ is uniformly bounded so that the oscillation rate remains bounded as $R\to \infty$. Combining this with an estimate for $\partial_t\bfu$ and using the equation along the approximating sequence, we eventually control the average of $\bfdelta$ and conclude existence of time-periodic weak solutions for the original problem in the whole $\Omega$.

All the above will be accomplished through several intermediate steps. From now on it will be tacitly understood that $\bfV(t)$ is $T$-periodic.


\subsection{The Initial-Boundary Value Problem}

In this subsection, we shall study  the following initial-boundary value problem associated to a mollified version of \eqref{3.1}:
\be
\tag{I-BVP}
\left\{\ba{l}\medskip\left.\ba{r}\medskip
\partial_t\bfu+\lambda(\bfu_\eta-{\bfgamma}+\bfU-\bfV)\cdot\nabla\bfu+(\bfu-\bfgamma)\cdot\nabla\bfU=\Div\mathbb T(\bfu,{p})+\bff\\
\Div\bfu=0\ea\right\}\ \ \mbox{in $\Omega_R\times(0,\infty)$}\,,\\ \medskip
\ \ \bfu(x,t)={\bfgamma}(t)\,, \ \mbox{ $(x,t)\in\partial\Omega\times(0,\infty)$}\,;\ \\ \medskip
\ \ \bfu(x,t)=\0\,, \ \mbox{ $(x,t)\in\partial B_R\times(0,\infty)$} \,;\ \\ \medskip
\ \ \bfu(x,0)=\bfu_0\,,\ x\in\Omega_R\,,\\ \medskip
\left.\ba{r}
\dot{\bfgamma}+\mathbb A\cdot\bfdelta+\varpi\Int{\partial\Omega}{} \mathbb T(\bfu,p)\cdot\bfn=\bfF\\ \medskip 
\dot{\bfdelta}=\bfgamma\ea\right\}\ \ \mbox{in $(0,\infty)$}\,,\\
\bfgamma(0)=\bfgamma_{0}\,,\ \ {\bfdelta}(0)=\bfdelta_0\,,
\ea\right.
\eeq{3.3}
where $\bfu_\eta$ is the (Friederichs) mollifier of $\bfu$ as defined in \eqref{regu}.
We will prove existence and uniqueness of solutions to \eqref{3.3} in a suitable functional class. To this end, we begin to transform it in an appropriate ``weak" form. Testing \eqref{3.3}$_1$ by $\bfpsi\in\cald^{1,2}(B_R)$, integrating by parts and using \eqref{3.3}$_{2-8}$
, we deduce
\be
\begin{split}
\langle\bfu(t)-\bfu_0,\bfpsi\rangle  
= &  -\Int0t\Big[\lambda\left((\bfu_\eta-{\bfgamma}+\bfU-\bfV)\cdot\nabla\bfu+(\bfu-\bfgamma)\cdot\nabla\bfU,\bfpsi\right)\\
&   +\left.2(\mathbb D(\bfu),\mathbb D ({\bfpsi}))+\varpi^{-1}{\hat{\bfpsi}}\cdot\mathbb A\cdot\bfdelta-(\bff,\bfpsi)-\bfF\cdot\hat{\bfpsi}\right]{\rm d}s\,,\\
\bfdelta(t) - \bfdelta_0 = &  \Int0t\bfgamma(s){\rm d}s\,,
\end{split}
\eeq{wefo}
where we recall that $(\cdot,\cdot)\equiv(\cdot,\cdot)_{\Omega_R}$ is the $L^2(\Omega_R)$-scalar product whereas 
$$
\langle \bfphi,\bfpsi\rangle=
\varpi^{-1}\,\hat{\bfphi}\cdot\hat{\bfpsi}+(\bfphi,\bfpsi)_{\Omega_R}.$$
\Bd The triple $(\bfu,\bfgamma,\bfdelta)$ is a \emph{weak solution} to \eqref{3.3} if, for all $t>0$:
\begin{itemize} 
\item [{\rm (i)}] $\bfu\in C_w([0,t]; L^2(\Omega_R))\cap L^2(0,t;\cald^{1,2}(\Omega_R))$\,, with $\bfu(x,t)|_{\partial\Omega}=\bfgamma(t)$\,,  \ \ $\bfgamma\in C([0,t];\real^3)$\,;\ 
\item [{\rm (ii)}] $\bfdelta\in C^{1}([0,t];\real^3)$ 
\,; 
\item [{\rm (iii)}] $(\bfu,\bfgamma,\bfdelta)$ satisfies \eqref{wefo} for all $\bfpsi\in\cald^{1,2}(B_R)$.
\end{itemize}
\EDD{4.1}
\Br
Taking into account classical properties of mollifiers and \remref{2.1}, it is easy to check that the integral in \eqref{wefo} is well defined for a weak solution. 
\ER{4.1}
\Br
With the help of the decomposition \eqref{Helm1},  one can show by classical arguments 
that if $(\bfu,\bfgamma,\bfdelta)$ satisfies \eqref{wefo} and is sufficiently regular, then there exists a suitable pressure field $p=p(x,t)$ such that $(\bfu,p,\bfgamma,\bfdelta)$ is a solution to \eqref{3.3}.
\ER{4.2}
 
The next lemma deals with the well-posedness of \eqref{3.3}.
For a weak solution $(\bfu,\bfgamma,\bfdelta)$ of \eqref{3.3}, we define the energy by
\be
E(t)=\half \left[\|\bfu(t)\|_{2,\Omega_R}^2+\varpi^{-1}\big(|\bfgamma(t)|^2+\bfdelta(t)\cdot\mathbb A\cdot\bfdelta(t)\big)\right].
\eeq{enE}
Observe  that, since $\bfgamma\equiv\hat{\bfu}$, the functional $E^{\frac12}$ defines a norm on $\calh(B_R)\times\real^3$. 

\Bl Let $\bfV\in W^{1,2}(0,T;\real^3)$, and let
\be
\eta_0<\exp\big(-8c_0\,\lambda \big)\,,\footnote{See Footnote (2).}
\eeq{eta}
with $\eta_0$ as in \lemmref{3.2}, $\eta\in(0,\eta_0)$ and $c_0$ as in \lemmref{1.1}--{\rm (vii)}.
Then, for any given $\bfu_0\in \calh(B_R)$, $\bfgamma_0\in\real^3$, $\bfdelta_0\in\real^3$, with the compatibility condition $\bfgamma_0:=\hat{\bfu_0}$, 
 there exists one and only one corresponding weak solution $(\bfu,\bfgamma,\bfdelta)$ to \eqref{3.3} such that 
\begin{itemize}
\item[\rm{(i)}] for all $\sigma>0$ and all $t>\sigma>0$,
\[
\bfu\in W^{1,2}(\sigma,t;L^2(\Omega_R))\cap L^2(\sigma,\tau;W^{2,2}(\Omega_R))\,,\ \ \bfgamma\in W^{1,2}(\sigma,t;\real^3)\,,\ \ \bfdelta\in W^{2,2}(\sigma,t;\real^3)\,;
\]
\item[\rm{(ii)}] for all $\sigma>0$ and all $t>\sigma>0$,
$\bfu\in C([\sigma,t];\cald^{1,2}(\Omega_R))$, and there exists $p\in L^2(\sigma,\tau;W^{1,2}(\Omega_R))$ such that $(\bfu,p,\bfgamma,\bfdelta)$ satisfies \eqref{3.3}$_1$ a.a. in $\Omega_R \times(0,t)$ and \eqref{3.3}$_{6,7}$ a.a. in $(0,t)$\,; 
\item[\rm{(iii)}] the initial conditions $(\bfu_0,\bfgamma_0,\bfdelta_0)$ are attained by $\bfgamma$ and $\bfdelta$ in the sense of pointwise continuity, and by $\bfu$ in the $L^2$-sense, i.e.
\be
\lim_{t\to 0}\|\bfu(t)-\bfu_0\|_{2,\Omega_R}=0\,;
\eeq{conIN} 
\item[\rm{(iv)}] for some $C_1,C_2$ independent of $R$ and $\eta$, some $C$ depending only on the data, $\eta$ and $t$,  the energy  estimates
\be	
\begin{split}
E(t)+C_1\Int0t(\|\nabla\bfu(s)\|^2_{2,\Omega_R}+|\bfgamma(s)|^2){\rm d}s & \le E(0)+C_2\Int0t(|\bfV|^2+|\dot{\bfV}|^2){\rm d}s\,,\ \ \mbox{for all $t>0$\,,}
\\
\Max{s\in[\sigma,t]}\|\nabla\bfu(s)\|_2 & \le  C\,,\ \ \mbox{for all $t>\sigma$,} 
\end{split}
\eeq{ESTi}
hold. 
\end{itemize}
Moreover, the solution depends continuously on the initial data in the norm $E^{\frac12}$.
\EL{3.3}

The proof of \lemmref{3.3} relies on Galerkin approximation. We will look for ``approximated" solutions to \eqref{wefo} of the form
\[
\bfu_N(x,t)= \sum_{k=1}^Nc_{kN}(t)\bfpsi_{k}(x)\,,\
\bfgamma_N(t)= \sum_{k=1}^Nc_{kN}(t)\hat{\bfpsi}_{k}\,,\ \bfdelta_N(t)\,, 
\]
where $\{\bfpsi_k\}$ is the special basis provided in the following lemma, whose proof is given in  \cite{GaSi1}. 
\Bl 
For any fixed $R>R_*$, the problem  
\be
\begin{array}{l}
\left. 
\begin{array}{l}
\hspace{-0.2cm} \displaystyle  - \nabla \cdot \mathbb T(\bfpsi,\phi)=\mu  \,\bfpsi \medskip \\ 
\hspace{-0.2cm} \displaystyle \Div   \bfpsi = 0 \medskip
\end{array}
\right\} \mbox{ in } \Omega_R \,,\\  \medskip
\displaystyle \bfpsi  = \hat{\bfpsi} \ \   \mbox{ in } \Omega_0\,, \ \  
\displaystyle \bfpsi= 0\ \   \mbox{ at } \partial B_R \,, \\ 
\displaystyle \mu\,  \hat{\bfpsi}=\varpi\int_{\partial\Omega} \mathbb T(\bfpsi,\phi)\cdot \bfn \,,
\end{array}
\eeq{eipr}
admits a denumerable number of positive eigenvalues $\{\mu_{i}\}$ clustering at infinity, 
and  corresponding eigenfunctions $\{\bfpsi_{i}\} \subset {\mathcal D}^{1,2}(B_R)\cap W^{2,2}(\Omega_R)$    forming an orthonormal basis of ${\mathcal{H}}(B_R)$ that is also orthogonal in $\cald^{1,2}(\Omega_R)$. Furthermore, the correspondent ``pressure" fields satisfy $\phi_{i}\in W^{1,2}(\Omega_R)$, $i\in\nat$.
\EL{Base}
\par

We also need the following approximation result that will be used in the proof of \propref{4.1}.
\Bl
Let $R>R_*$ be given, and let $\{\bfpsi_{k}\}_{k \in \mathbb{N}}$ be the basis given in 
\lemmref{Base}. Then, for any $\bfphi\in C_\sharp(B_R)$  and any $\varepsilon >0$ there
is $N=N(\bfphi,\varepsilon )\in {\mathbb N}$, and corresponding $T$-periodic functions $\ds r_{k}\in
C^{1}(\real),$ $k=1,...,N,$ such that
$$
\Max{t\in [0,T]} \left\{ \| (\bfphi _{N}-\bfphi)(.,t) \|_{1,2} + |(\bfphi_{1N}-\bfphi_1 )(t)| 
+ \| \partial_t (\bfphi _{N} - \bfphi )(.,t) \|_{2} + |\partial_t(\bfphi_{1N}-\bfphi_1 )(t)|
\right\} < \varepsilon 
$$
with $\ds \bfphi _{N}(x,t)=\sum_{k=1}^{N}r_{k}(t)\bfpsi_{k}(x)$, and $\ds \bfphi _{1N}(t)=\sum_{k=1}^{N}r_{k}(t)\hat{\bfpsi}_{k}$. 
\EL{2.3}
The proof is obtained arguing as in \cite[Lemma 3.1]{GaSi2} and therefore omitted.

\medbreak

As already said, \lemmref{3.3} will be proved with the classical Galerkin procedure. To start the process, we search for an approximated solution $(\bfu_N(x,t), \bfgamma_N(t), \bfdelta_N(t))$ to \eqref{wefo} of the form
\be
\begin{split}
\bfu_N(x,t)= &\sum_{k=1}^Nc_{kN}(t)\bfpsi_{k}(x)\,,\\ 
\bfgamma_N(t)= &\sum_{k=1}^Nc_{kN}(t)\hat{\bfpsi}_{k}\,, 
\end{split}
\eeq{gal0}
where $\{\bfpsi_k\}$ is the basis introduced in \lemmref{Base} and the vector functions $\bfc_N(t):=\{c_{1N}(t),\ldots c_{NN}(t)\}$ and $\bfdelta_N(t)$ satisfy the following system of  equations 
\be
\ba{ll}\medskip
\ode{}t\langle\bfu_N,\bfpsi_i\rangle+\lambda\left(((\bfu_N)_\eta-{\bfgamma_N}+\bfU-\bfV)\cdot\nabla\bfu_N+(\bfu_N-\bfgamma_N)\cdot\nabla\bfU,\bfpsi_{i}\right)\\ \medskip\hspace*{1.5cm}=-2(\mathbb D(\bfu_N),\mathbb D ({\bfpsi}_{i}))-\varpi^{-1}\hat{\bfpsi}_i\cdot\mathbb A\cdot\bfdelta_N+(\bff,\bfpsi_i)+\bfF\cdot\hat{\bfpsi}_i\,,
\\
\hspace*{.9cm}\dot{\bfdelta}_N=\bfgamma_N\,,
\ea
\eeq{gal}
 $i=1,\ldots,N$. 
Here 
$(\cdot,\cdot)\equiv(\cdot,\cdot)_{\Omega_R}$ and $\langle\cdot,\cdot\rangle\equiv\langle\cdot,\cdot\rangle_{B_R}$. This yields a system of first order differential equations in normal form in the unknowns $\bfc_N,\bfchi_N$. Indeed, since we have the orthogonality conditions
\be
\langle\bfpsi_{i},\bfpsi_{j}\rangle=\delta_{ij}\,,
\eeq{ortg}
plugging the ansatz \eqref{gal0} into \eqref{gal} entails
\be
\begin{split}
\dot{c}_{iN}= &\ {\sf F}_i(\bfc_N,\bfdelta_N)\,,\ \ i=1,\ldots,N\,;\,\\ 
\dot{\bfdelta}_N= &\ \bfgamma_N\,, 
\end{split}
\eeq{ODE}
where
\be\ba{ll}\medskip
{\sf F}_i:=-\!\Sum{k=1}Nc_{kN}\left[\lambda\big((\!\bfU\!-\!\bfV)\cdot\nabla\bfpsi_k+(\bfpsi_k\!-\!\hat{\bfpsi}_k)\cdot\nabla\bfU,\bfpsi_{i}\big)\!+\!2\big(\mathbb D(\bfpsi_{k}),\mathbb D(\bfpsi_{i})\big)\right]\\
\hspace*{1.7cm}-\Sum{k=1}N\Sum{l=1}Nc_{kN}c_{lN}(((\bfpsi_{k})_\eta-\hat{\bfpsi}_{k})\cdot\nabla\bfpsi_l,\bfpsi_i)-\Frac{1}{\varpi}\hat{\bfpsi}_{i}\cdot\mathbb A\cdot\bfdelta_N+(\bff,\bfpsi_i)+\bfF\cdot\hat{\bfpsi}_i\,.
\ea
\eeq{ODE1}
The initial conditions $(\bfu_0,\bfgamma_0,\bfdelta_0)$ at the level of the coefficients read
\be
c_{iN}(0)= \langle\bfu_0,\bfpsi_{i}\rangle =(\bfu_0,\bfpsi_{i})+\varpi^{-1}\bfgamma_0\cdot\hat{\bfpsi}_{i}\,,\ \ \bfdelta_N(0)=\bfdelta_0\,.
\eeq{inco}
Since $\{\bfpsi_{i}\}$ is an orthonormal basis of ${\mathcal{H}}(B_R)$, multiplying the first identity in \eqref{inco} by $\bfpsi_{i}$ and summing over the index $i$ from $1$ to $N$ delivers a bound on the initial conditions $(\bfu_N(0),\bfgamma_N(0),\bfdelta_N(0))$
\be
\|\bfu_N(0)\|_{2}^2+\varpi^{-1}|\bfgamma_N(0)|^2\le\|\bfu_0\|_{2}^2+\varpi^{-1}|\bfgamma_0|^2\,.
\eeq{incou}
Now that we have settled the starting point of the argument, we turn to the proof of \lemmref{3.3}. 

\medbreak

\noindent{\em Proof of \lemmref{3.3}.}
We begin to derive three basic energy estimates for the approximated solution $\bfu_N$ which are the approximated forms of estimates that would be obtained formally for $\bfu$ by choosing abusively $\bfu$, $-t\,\Div\mathbb T(\bfu,p)$ and $t\,\partial_t\bfu$ as test functions in \eqref{wefo}. Once we have these estimates at hand, we can actually let $N$ go to infinity and prove the assertions \rm{(i)--(iv)} for the limit $\bfu$ of the sequence $\{\bfu_N\}$. Then it will remain to prove the continuous dependence on the initial data implying at once the uniqueness of the solution $\bfu$. 

\medbreak

\noindent{\bf Energy estimates.}
{\it First estimate:}
to mimic the choice of $\bfpsi=\bfu$ in \eqref{wefo}, we multiply both sides of \eqref{gal}$_1$ by $c_{iN}$, sum over $i$ and integrate by parts over $\Omega_R$. Using \eqref{gal0} along with \lemmref{3.2}--(c), we show
\be
\begin{split}
\half\ode{}t\left[\|\bfu_N\|_2^2+\varpi^{-1}(|\bfgamma_N|^2+\bfdelta_N\cdot\mathbb A\cdot\bfdelta_N)\right] & +2\|\mathbb D(\bfu_N)\|_2^2 \\ & =\lambda\,((\bfu_N-\bfgamma_N)\cdot\nabla\bfu_N,\bfU)+(\bff,\bfu_N)+\bfF\cdot\bfgamma_N\,. 
\end{split}
\eeq{3.12}
Since $\eta_0$ satisfies \eqref{eta} and $c_0$ does not depend on $\varepsilon$, we can choose $\varepsilon$ in the construction of $\bfU$ such that 
\be
\lambda|((\bfu_N-\bfgamma_N)\cdot\nabla\bfu_N,\bfU)|\le \lambda c_0\varepsilon \|\nabla\bfu_N\|_2^2 \le  \frac{1}{4}\,\|\nabla\bfu_N\|_2^2\,,
\eeq{3.13}
whereas, from \eqref{1.8} and \eqref{1.10} it follows
\be
2\|\mathbb D(\bfu_N)\|_2^2-\mbox{$\frac14$}\|\nabla\bfu_N\|_2^2\ge \half\|\nabla\bfu_N\|_2^2+\kappa\,|\bfgamma_N|^2\,.
\eeq{3.14}
for some universal $\kappa>0$. 

We next estimate the last terms in \eqref{3.12}. Recalling that  $\bff=(\bfV-\bfU)\cdot\nabla\bfU+\Delta\bfU-\partial_t\bfU$ and 
$\bfF=-\varpi\,{\rm vol}\,(\mathscr B)\dot{\bfV}$, the properties of $\bfU$ proved in \lemmref{1.1}, together with H\"older inequality, \eqref{1.6} and Cauchy-Schwarz inequality implies
$$
|(\bff,\bfu_N)|+|\bfF\cdot\bfgamma_N|\le c\,(|\bfV(t)|^2+|\dot{\bfV}(t)|^2)+\mbox{$\frac14$}\|\nabla\bfu_N\|_2^2+\mbox{$\frac12$}\kappa|\bfgamma_N|^2\,.
$$ 
Thus, employing in \eqref{3.12} the latter together with \eqref{3.13}, \eqref{3.14}, we establish that
\be
\half\ode{}t\left[\|\bfu_N\|_2^2+\varpi^{-1}(|\bfgamma_N|^2+\bfdelta_N\cdot\mathbb A\cdot\bfdelta_N)\right]+\mbox{$\frac14$}\|\nabla\bfu_N\|_2^2+\half\kappa|\bfgamma_N|^2\le c_2\,(|\bfV(t)|^2+|\dot{\bfV}(t)|^2)\,.
\eeq{3.15}
Denoting by $E_N(t)$ the energy of the approximated solution, i.e.
\be
E_N(t):=\half\left[\|\bfu_N(t)\|_2^2+\varpi^{-1}(|\bfgamma_N(t)|^2+\bfdelta_N(t)\cdot\mathbb A\cdot\bfdelta_N(t))\right]\,,
\eeq{3.16}
and using the bound \eqref{incou} on the initial conditions, we infer 
\be
\begin{split}
E_N(t) & \ +c_3\Int0t(\|\nabla\bfu_N(s)\|^2  +|\bfgamma_N(s)|^2){\rm d}s  \le E_N(0)+c_2\Int0t\left(|\bfV(s)|^2+|\dot\bfV(s)|^2\right){\rm d}s\\
& \ \le \half(\|\bfu_0\|^2+\varpi^{-1}\big(|\bfgamma_0|^2+\bfdelta_0\cdot\mathbb A\cdot\bfdelta_0)\big)+c_2\Int0t\left(|\bfV(s)|^2+|\dot\bfV(s)|^2\right){\rm d}s\,,
\end{split}
\eeq{3.17}
where $c_2$ and $c_3$ are independent of $R$ and $\eta$.

As a consequence of this energy bound,  we deduce that whatever $\tau>0$, there exists $C(\tau)>0$, independent of $R$ and $\eta$, such that 
$$|\bfc_N(t)|+|\bfdelta_N(t)|\le C(\tau), \quad t\in[0,\tau].$$ In particular, this bound in turn implies that the initial-value problem \eqref{ODE}--\eqref{inco} has a unique global solution (i.e. defined for all $t>0$). 

\medbreak

\noindent{\it Second estimate : }
to mimic the formal choice $\bfpsi=-t\,\Div\mathbb T(\bfu,p)$ in \eqref{wefo}, we next multiply both sides of \eqref{gal}$_1$ by $t\,\mu_{i}c_{iN}$, $t>0$, and sum over $i$. Integrating by parts over $\Omega_R$ and employing \eqref{eipr}, we show
\be\begin{split}
\half\ode{}t\Big( &t\|\nabla \bfu_N\|_2^2\Big)+ t\|\Div\mathbb T(\bfu_N,p_N)\|_2^2+\varpi\,t\,|\bfS_N|^2  -\half\|\nabla\bfu_N\|_2^2\\ \medskip
= &\ \lambda\,t\,\left((\bfu_N)_\eta-{\bfgamma_N}+\bfU-\bfV)\cdot\nabla\bfu_N+(\bfu_N-\bfgamma_N)\cdot\nabla\bfU,\Div\mathbb T(\bfu_N,p_N)\right)
+t\,\bfS_N\cdot\mathbb A\cdot\bfdelta_N\\
&\ +t\,(\bff,\Div\mathbb T(\bfu_N,p_N))+t\,\varpi\,(\bfF,\bfS_N)
 \end{split}
\eeq{3.18}
where
$$
\bfS_N:=\int_{\partial\Omega}\mathbb T(\bfu_N,p_N)\cdot\bfn\,,\ \ p_N:=\sum_{k=1}^Nc_{kN}\phi_{k}\,,
$$
and $\phi_{k}$ is the ``pressure" field associated to $\bfpsi_{k}$. 
We estimate piece by piece the right-hand side of \eqref{3.18}. Since
$$
\|(\bfu_N)_\eta\|_\infty\le c_\eta\|\bfu_N\|_2,
$$
for some $c_\eta>0$, we have 
\be
\lambda\,|\left((\bfu_N)_\eta\cdot\nabla\bfu_N,\Div\mathbb T(\bfu_N,p_N)\right)|\le c_\eta\|\bfu_N\|_2^2\|\nabla\bfu_N\|_2^2+\mbox{$\frac18$}\|\Div\mathbb T(\bfu_N,p_N)\|_2^2\,.
\eeq{3.20}
Using H\"older inequality,  the properties of $\bfU$, \eqref{1.9} and \eqref{1.10}, we also get
\be\ba{rl}\medskip
\lambda|\left((\bfu_N-\bfgamma_N)\cdot\nabla\bfU,\Div\mathbb T(\bfu_N,p_N)\right)| \le & c\,|\bfV(t)|^2\|\nabla\bfu_N\|_2^2+\frac18\|\Div\mathbb T(\bfu_N,p_N)\|_2^2\\  \medskip
\lambda|\left(-{\bfgamma_N}+\bfU-\bfV)\cdot\nabla\bfu_N,\Div\mathbb T(\bfu_N,p_N)\right)|\le & c\,(\|\nabla\bfu_N\|_2^2+|\bfV(t)|^2)\|\nabla\bfu_N\|_2^2\\ \medskip 
& +\frac18\|\Div\mathbb T(\bfu_N,p_N)\|_2^2\,. 
\ea
\eeq{3.19}
We now set
$$
g(t):=1+\|\bfu_N(t)\|_{1,2}^2+|\bfV(t)|^2\,,$$ 
and 
$$
G(t):= \|\nabla\bfu_N(t)\|_2^2+t|\bfdelta_N(t)|^2+t\,\|\bff(t)\|_2^2+t\,|\bfF(t)|^2
$$
and we observe that from the energy estimate \eqref{3.17}, recalling also \eqref{1.7}, we obtain for an arbitrary $\tau>0$, 
$$
\int_0^\tau(g(t)+G(t)){\rm d}t\le F(\tau),
$$
where $F$ depends only on the norm of the initial data, the $W^{1,2}$-norm of $\bfV(t)$, $\tau$, and $\eta$.  
Finally, we make use of the classical estimate for the Stokes problem
\be
\|\bfu_N\|_{2,2}\le c_5\,(\|\Div\mathbb T(\bfu_N,p_N)\|_2+|\hat{\bfu_N}|)\le c_6\,(\|\Div\mathbb T(\bfu_N,p_N)\|_2+\|\nabla\bfu_N\|_2)\,,
\eeq{3.21}
where, in the last inequality, we have used \eqref{1.10}, and the constant $c_6$ depends on $R$.
Then, using one more time Cauchy-Schwarz inequality, we conclude from \eqref{3.18}--\eqref{3.21} that
\be
\ode{}t(t\|\nabla\bfu_N\|^2_2)+c_7\,t\|\bfu_N\|_2^2\le c_8\,g(t)(t\|\nabla\bfu_N\|_2^2)+c_9\,G(t)\,, 
\eeq{3.22}
the constants depend on $R$.

As a result, using Gronwall's lemma in \eqref{3.22} entails
\be
\sup_{t\in [0,\tau]}\left(t\,\|\nabla\bfu_N(t)\|_2\right)+\int_0^\tau(t\|\bfu_N(t)\|_{2,2}^2){\rm d}t\le H_1(\tau)\,,\ \ \mbox{for all $\tau>0$}\,,
\eeq{3.23}
where $H_1$ has the same property as $F$ and, in addition, depends also on $R$.

\medbreak

\noindent{\it Third estimate: }
we finally mimic the formal choice $\bfpsi=t\,\partial_t\bfu$ in \eqref{wefo} multiplying both sides of \eqref{gal}$_1$ by $t\,\dot{c}_{iN}$, summing over $i$,  and integrating by parts over $\Omega_R$ as necessary. Taking again into account \eqref{gal0}, \eqref{1.8} and \remref{2.1}, we show
\be\ba{ll}\medskip
\half\ode{}t(t\|\nabla\bfu_N\|_2^2)+t\|\partial_t\bfu_N\|_2^2+\varpi^{-1}t|\dot{\bfgamma}_N|^2-\frac12\|\nabla\bfu_N\|_2^2\\ \medbreak \quad =\lambda\,t\,\left((\bfu_N)_\eta - {\bfgamma_N}+\bfU\!-\!\bfV)\cdot\nabla\bfu_N+(\bfu_N-\bfgamma_N)\cdot\nabla\bfU,\partial_t\bfu_N\right)\\ \medbreak
\quad -\frac{t}{\varpi}\dot{\bfgamma}_N\cdot\mathbb A\cdot\bfdelta_N+t(\bff,\partial_t\bfu_N)+t\bfF\cdot\dot{\bfgamma}_N\,.
\ea
\eeq{3.24}
We may now proceed to estimate the right-hand side of \eqref{3.24} exactly as as we did in  \eqref{3.19}, \eqref{3.20}, with $\partial_t\bfu_N$ in place of $\Div\mathbb T$, to show that \eqref{3.24} implies the following further bound
\be 
\int_0^\tau t(\|\partial_t\bfu_N(t)\|_{2}^2+|\dot{\bfgamma}_N(t)|^2){\rm d}t\le H_2(\tau)\,,\ \ \mbox{for all $\tau>0$}\,,
\eeq{3.25}
where $H_2$ has the same property as $H_1$. 

\medbreak
\noindent{\bf Convergence of the sequence \mbox{\boldmath$\{\bfu_N,\bfgamma_N,\bfdelta_N\}$}.} 
We warn the reader that, along this step of the proof, convergence is understood up to the choice of an ad-hoc subsequence if needed, and all extracted subsequences will still be denoted by $\{\bfu_N,\bfgamma_N,\bfdelta_N\}$ for simplicity. 
Integrating both sides of \eqref{gal}$_1$ over the interval $[t_1,t_2]$, $0\le t_1<t_2$,  and taking into account \eqref{ortg}, we get
\be\ba{ll}\medskip
\langle\bfu_N(t_2)-\bfu_N(t_1),\bfpsi_i\rangle=
-\Int{t_1}{t_2}\big\{\lambda\left(((\bfu_N)_\eta-{\bfgamma_N}+\bfU-\bfV)\cdot\nabla\bfu_N+(\bfu_N-\bfgamma_N)\cdot\nabla\bfU,\bfpsi_{i}\right)\\ \medskip\hspace*{5.5cm}+2(\mathbb D(\bfu_N),\mathbb D ({\bfpsi}_{i}))+\varpi^{-1}\mathbb A\cdot\bfdelta_N\cdot{\hat{\bfpsi}}_{i}-(\bff,\bfpsi_i)-\bfF\cdot\hat{\bfpsi}_i\big\}\,{\rm d}t\,.
\ea
\eeq{3.26}
Employing \eqref{3.26} and the uniform bound \eqref{3.17}, by means of a classical argument \cite[Section 3]{Gaintro}, we can show the existence of 
\be
\bfu\in L^\infty(0,t;\calh(B_R))\ \mbox{for all $t>0$}\,, 
\eeq{Linf}
such that for all $\bfpsi\in\cald^{1,2}(B_R)$, and setting $\bfgamma=\hat{\bfu}$, 
\be\ba{ll}\medskip
\Lim{N\to\infty}\big[(\bfu_N(t),\bfpsi)+\varpi^{-1}\bfgamma_N(t)\cdot\hat{\bfpsi}\big]=(\bfu(t),\bfpsi)+\varpi^{-1}\bfgamma(t)\cdot\hat{\bfpsi}\,,\ea
\eeq{3.27}
uniformly in $t\in [0,\tau]$, for all $\tau>0$.
This implies, in particular,
\be
\Lim{N\to\infty}\bfgamma_N(t)=\bfgamma(t)\,,\ \ \mbox{uniformly in $t\in [0,\tau]$}\,,
\eeq{3.28}
which  furnishes 
\be
\bfgamma\in C([0,\tau];\real^3)\,,
\eeq{gaco} 
and also, by \eqref{ODE}$_2$,
\be
\bfdelta_N\to\bfdelta\ \ \mbox{in $C^1([0,\tau];\real^3)$}\,,\ \ \dot{\bfdelta}(t)=\bfgamma(t)\,,\ \ t\in [0,\tau]\,.
\eeq{3.29}
From the uniformity in $N$ of the estimates \eqref{3.17}, \eqref{3.23} and \eqref{3.25}, it is routine to deduce that 
\be
\bfu_N\to\bfu\ \left\{\ba{ll}\medskip \mbox{weakly in $L^2(0,t;W^{1,2}(\Omega_R))$}\,,\\ \mbox{weak$^*$ in $L^\infty(0,t;L^2(\Omega_R))$\,,}\ea\right. 
\eeq{3.30}
for all $t>0$ and so, combining \eqref{3.27} with \eqref{3.30}$_1$, we also deduce (we refer to \cite[Section 5]{Gaintro} for details)
\be
\bfu_N\to\bfu\,, \ \mbox{strongly in $L^2(0,t;L^{2}(\Omega_R))$}\,,\ \ \mbox{for all $t>0$}.
\eeq{3.31}
In particular, from \eqref{3.30} \eqref{3.17}, \eqref{3.28} and \eqref{3.29} we infer the validity of \eqref{ESTi}$_1$.  
Likewise, from \eqref{3.23} and \eqref{3.25} it follows that
\be
\bfu_N\to \bfu \ \ \mbox{weakly in $W^{1,2}(\sigma,\tau;\call^{2}(B_R))$ and $L^2(\sigma,\tau;W^{2,2}(\Omega_R))$, for all $\tau>\sigma>0$\,.}
\eeq{3.32}
The latter implies, on the one hand, by a well-known interpolation theorem,  $\bfu\in C([\sigma,\tau];W^{1,2}(\Omega_R))$, and, on the other hand, in combination with \eqref{3.23}, the estimate \eqref{ESTi}$_2$. 

We next pass to the limit $N\to\infty$ in \eqref{3.26}. Using \eqref{3.27}, \eqref{3.30} and \eqref{3.31} together  with the properties of the base $\{\bfpsi_i\}$ (see \lemmref{Base}) and classical arguments from \cite[Section 5]{Gaintro}, we can show that \eqref{3.26} continues to hold with $(\bfu_N,\bfgamma_N,\bfdelta_N)$ replaced by $(\bfu,\bfgamma,\bfdelta)$ and $\bfpsi_i$ replaced by arbitrary $\bfpsi\in\cald^{1,2}(B_R)$. 

All in all, $(\bfu,\bfgamma,\bfdelta)$ satisfies \eqref{wefo}, for all $t>0$ and in view of \eqref{gaco} and \eqref{Linf}, this implies also 
\be\bfu\in C_w(0,t;L^2(\Omega_R))\,.
\eeq{uco}
In order to complete the existence part of the proof, it remains to consider the sense in which the initial conditions $(\bfu_0,\bfgamma_0,\bfdelta_0)$ are satisfied and in fact it remains only to show \eqref{conIN}. The energy estimate \eqref{ESTi} and the continuity of $(\bfgamma,\bfdelta)$ at $t=0$ imply 
$$
\limsup_{t\to 0^+}\|\bfu(t)\|_2\le\|\bfu_0\|_2\,,
$$
which entails \eqref{conIN} as immediate consequence of the weak continuity \eqref{uco}.

\medbreak
\noindent{\bf Continuous dependence on the initial data.}
Let $(\bfu_i,\bfgamma_i,\bfdelta_i)$, $i=1,2$, be two weak solutions corresponding to the same \textcolor{black}{$\bfV$}. Setting $\bfu=\bfu_1-\bfu_2$, $\bfgamma=\bfgamma_1-\bfgamma_2$, $\bfdelta=\bfdelta_1-\bfdelta_2$ we have for arbitrary $\tau>\sigma>0$
\be
\left\{\ba{l}\medskip\left.\ba{r}\medskip
\partial_t\bfu\!+\!\lambda((\bfu+\bfu_2)_\eta\!-\!{\bfgamma}+\bfU\!-\!\bfV)\cdot\nabla\bfu\!+\!(\bfu\!-\!\bfgamma)\cdot\nabla\bfU\!+\!\bfu_\eta\cdot\nabla\bfu_2=\Div\mathbb T(\bfu,{p})\\
\Div\bfu=0\ea\right\}\ \, \mbox{ in $\Omega_R\times[\sigma,\tau]$}\,,\\ \medskip
\ \ \bfu(x,t)={\bfgamma}(t)\,, \ \mbox{ $(x,t)\in\partial\Omega\times[\sigma,\tau]$}\,,\ \  \\ \medskip
\ \ \bfu(x,t)=\0\,, \ \mbox{a.a. $(x,t)\in\partial B_R\times[\sigma,\tau]$}\,, \\ \medskip
\ \ \bfu(x,0)=\bfu_1(x,0)-\bfu_2(x,0)\,,\ x\in\Omega_R\,,\\ \medskip
\left.\ba{r}
\dot{\bfgamma}+\mathbb A\cdot\bfdelta+\varpi\Int{\partial\Omega}{} \mathbb T(\bfv,p)\cdot\bfn=\0\\ \medskip 
\dot{\bfdelta}=\bfgamma\ea\right\}\ \ \mbox{in $(0,\infty)$}\,,\\
\bfgamma(0)=\bfgamma_1(0)-\bfgamma_2(0)\,,\ \ {\bfdelta}(0)=\bfdelta_1(0)-\bfdelta_2(0)\,,
\ea\right.
\eeq{4.40}
for some $p\in L^2(\sigma,\tau;W^{1,2}(\Omega_R))$. Testing \eqref{4.40}$_1$ by $\bfu$, integrating by parts over $\Omega_R$ and taking into account  \eqref{4.40}$_{2-5}$ we get
\be
\half\ode{E}t +2\|\mathbb D(\bfu)\|_2^2=\lambda\,[((\bfu-\bfgamma)\cdot\nabla\bfu,\bfU)-(\bfu_\eta\cdot\nabla\bfu_2,\bfu)]\,,\ t\in[\sigma,\tau]\,,
\eeq{4.41}
where $E$ is given in \eqref{enE}.
Arguing as in \eqref{3.13}, \eqref{3.14} from \eqref{4.41} we infer, in particular,
$$
\ode{E}t\le -\lambda\,(\bfu_\eta\cdot\nabla\bfu_2,\bfu)\,,\ t\in[\sigma,\tau].
$$
By the property of mollifier, we have $\|\bfu_\eta\|_\infty\le c_\eta\|\bfu\|_2$, so that
$$
|(\bfu_\eta\cdot\nabla\bfu_2,\bfu)|\le c_\eta\|\nabla\bfu_2\|_2\|\bfu\|_2^2\le2 c_\eta \|\nabla\bfu_2\|_2 E\,.
$$
Combining the last two displayed equations and using Gronwall's lemma, we get
$$
E(\tau)\le E(\sigma)\exp(2c_\eta\int_0^\tau\|\nabla\bfu_2(s)\|_2{\rm d}s)\,.
$$
If we let $\sigma\to0$ and use the properties (i) and (ii) of a weak solution along with \eqref{conIN}, we get
$$\lim_{\sigma\to 0}E(\sigma)= E(0)
$$
yielding the claimed continuous dependence property and, in particular, uniqueness for the Cauchy problem.
\par\hfill$\square$\par

We end this subsection on the initial-boundary value problem with the following important estimate that will be used later on in \propref{4.1}.
\Bl Let $\bfV\in W^{1,2}(0,T;\real^3)$. Then, every weak solution to \eqref{3.3} satisfies for all $t>0$ and all $R>R_*$
$$
\ode{\bfu}{t}\in L^1(0,t;\cald^{-1,2}(B_R))\,.
$$
Moreover, there exists a constant $C=C(R)$ independent of $\eta$ such that
\be
\left\|\ode{\bfu}{t}\right\|_{L^1(0,T;\cald^{-1,2}(B_{R}))}\le C\,\int_0^t\big(|\bfV(s)|+|\dot{\bfV}(s)|+|\bfdelta(s)|^2+\|\nabla \bfu(s)\|_{2}+\|\nabla \bfu(s)\|_{2}^2){\rm d}s\,.
\eeq{4.53} 
\EL{4.4}
{\em Proof.} From \eqref{wefo}$_1$ we deduce, for a.a. $t>0$ and all $\bfpsi\in \cald^{1,2}(B_{R})$
\be
\ode{}t\langle \bfu(t),\bfpsi\rangle=G_{\mbox{\tiny $\bfpsi$}}(t)\,,
\eeq{4.54}
where
$$
G_{\mbox{\tiny $\bfpsi$}}(t)=-\lambda\left((\bfu_\eta-{\bfgamma}+\bfU-\bfV)\cdot\nabla\bfu+(\bfu-\bfgamma)\cdot\nabla\bfU,\bfpsi\right)-2(\mathbb D(\bfu),\mathbb D ({\bfpsi}))-\varpi^{-1}{\hat{\bfpsi}}\cdot\mathbb A\cdot\bfdelta+(\bff,\bfpsi)+\bfF\cdot\hat{\bfpsi}
\,.
$$
Employing several times H\"older inequality along with \eqref{1.8}--\eqref{1.10} and classical properties of the mollifier, we show 
\be
|G_{\mbox{\tiny $\bfpsi$}}(t)|\le C(R) \,\big(|\bfV(t)|+|\dot{\bfV}(t)|+|\bfdelta(t)|^2+\|\nabla \bfu(t)\|_{2}+\|\nabla \bfu(t)\|_{2}^2)\|\mathbb D(\bfpsi)\|_{2} 
\,.
\eeq{4.55}
Thus,  from \eqref{4.54}--\eqref{4.55} and the assertions (i)-(ii) of \defref{4.1}, we deduce that there exists $\bfg\in L^1(0,T; \cald^{-1,2}(B_{R}))$ such that
$$
\ode{}t\bfu(t)= g(t)
$$
in the sense of distributions,  as well as the validity of \eqref{4.53}. 
\par\hfill$\square$\par

\subsection{A Key Lemma on the Total Dissipation}

As explained in the introduction, we will obtain the existence of a periodic solution as fixed point of the Poincar\'e map  ${\sf M}$ in the energy space. To prove that ${\sf M}$ maps some ball in the energy space into itself, we next adapt an argument originally introduced in the context of nonlinear wave equations \cite[p. 162]{HX},\cite{HX2}.  The basic idea, reformulated in our case, consists in perturbing the energy by adding the term $\zeta \varpi^{-1}\bfgamma\cdot\bfdelta$ for the solid, to recover some dissipation in $\bfdelta$. To be compatible with the fluid, we need to add another ad-hoc contribution in the energy functional. To this aim, we first construct (classically) a solenoidal extension of $\bfdelta=\bfdelta(t)$. Set
$$\ba{ll}\medskip
\bsfH(x,t):=x_3\delta_2(t)\bfe_1+x_1\delta_3(t)\bfe_2+x_2\delta_1(t)\bfe_3
\\
\bfH(x,t):=\curl(\psi(|x|)\bsfH(x,t))\equiv \psi(|x|)\bfdelta(t)-\bsfH(x,t)\times\nabla\psi(|x|)\,,
\ea$$
where $\psi$ is a smooth function that takes value $1$ in a neighborhood of $\Omega_0$ and vanishes for $|x|\ge{2R_*}$. Clearly, 
\be\ba{ll}\medskip
\bfH(x,t)=\bfdelta(t)\,,\, \ (x,t)\in\partial\Omega\times (0,\infty)\,;\ \\ \medskip 
\Div\bfH(x,t)=0\,,\,\ (x,t)\in\Omega_R\times(0,\infty)\,;\\ \medskip
\Sup{x\in\Omega_R}\left|\bfH(x,t)\right|\le c\,|\bfdelta(t)|\,,\\ \medskip
\Sup{x\in\Omega_R}\left|\partial_t\bfH(x,t)\right|\le c\,|\dot{\bfdelta}(t)|=c\,|\bfgamma(t)|\,.    
\ea
\eeq{deH}
We then introduce the map from $C^{1}([0,t];\real^3)$ to $L^2(0,t;\cald^{1,2}(\Omega_R))$ that associate the solenoidal extension $\bfH(t)$ we just built to $\bfdelta(t)$. We could write this field $\bfH_{\mbox{\footnotesize $\bfdelta$}}(t)$ but we just keep the notation $\bfH(t)$ since no confusion arises.
Let $(\bfu,\bfgamma,\bfdelta)$ be the weak solution of \eqref{3.3} determined in \lemmref{3.3}. For $\zeta\in (0,\infty)$, we define the following modified energy functional
$$
E_\zeta(\bfu,\bfgamma,\bfdelta):= E(\bfu,\bfgamma,\bfdelta)+ 2\zeta\,[(\bfu,\bfH)+\varpi^{-1}\bfgamma\cdot\bfdelta]\,,
$$
where $E$ is given in \eqref{enE}. With the help of \eqref{deH} and Cauchy-Schwarz inequality, it is easy to check that there  is $\zeta_1=\zeta_1(\varpi,R_*,\mathbb A)$ such that, if $\zeta\le\zeta_1$, then
\be 
\half E\le E_\zeta\le 2 E\,.
\eeq{4.44}
  As a consequence,  $E^{\frac12}$ and $E_\zeta^{\frac12}$ are equivalent norms in $\call^2(B_R)$ and the set 
$\{(\bfu,\bfdelta)\in \calh(B_R)\times\real^3:\, \ E_{\zeta}(\bfu,\hat{\bfu},\bfdelta)\le \rho^2\}$ is convex for $\zeta$ small. This fact is needed in the next subsection, in order  to obtain a fixed point of the Poincar\'e map via Schauder theorem. 
 
We next show that there exist $\zeta_0,\rho_0>0$,  depending on the parameters of the problem such that, for the norm $E_{\zeta_0}^{\frac12}$, the Poincar\'e map maps the ball of radius $\rho_0$ in $\call^2(B_R)\times\real^3\times\real^3$ into itself.
To this end, we set
$$
\mathcal V:=\int_0^T(|\bfV|^2+|\dot{\bfV}|^2){\rm d}t\,.
$$

\Bl Let $\bfV$ satisfy the assumption of   \lemmref{3.3} and let $(\bfu,\bfgamma,\bfdelta)$ be the corresponding unique solution to  \eqref{3.3} given in \lemmref{3.3}. 
There exist $\zeta_0,\rho_0>0$,  depending on $R$, $\mathcal V$, $\lambda$, $\mathbb A$ and $\varpi$, such that if $E_{\zeta_0}(0)\le \rho_0$, then $E_{\zeta_0}(T)\le \rho_0$ as well. In addition, the following estimate holds
\be
E_{\zeta_0}(T)+\half\zeta_0\varpi^{-1}\!\!\int_0^T\bfdelta(t)\cdot\mathbb A\cdot\bfdelta(t)\,{\rm d}t\le E_{\zeta_0}(0)+C_0
\mathcal V\,,
\eeq{Ez}
where $C_0$ is independent of $\eta$.
\EL{4.3}
{\em Proof.} For sufficiently small $\sigma>0$ consider \eqref{3.3} with  $t\in[\sigma,T]$. Testing \eqref{3.3}$_1$ by $\bfu$ and using exactly the same argument leading to \eqref{3.15} we obtain
\be
\ode{E}t+\half\|\nabla\bfu\|_2^2+\kappa|\bfgamma|^2\le C_1\,(|\bfV|^2+|\dot{\bfV}|^2)\,,
\eeq{4.45}
where, here and in the rest of the proof, by $C_i$, $i=1,\ldots$, and, later on, $\zeta_2$ we denote  generic positive constants depending, at most, on $\lambda,\mathbb A,$ $\varpi$ and the (spatial)  support of $\bfH$ but not on $\eta$. Furthermore, again from \eqref{3.3}$_1$, tested this time on $\bfH$, we get
\be\ba{ll}\medskip
\ode{}t[(\bfu,\bfH)+\varpi^{-1}\bfgamma\cdot\bfdelta]= (\bfu,\partial_t\bfH)+\varpi^{-1}|\bfgamma|^2-\lambda\left((\bfu_\eta-{\bfgamma}+\bfU-\bfV)\cdot\nabla\bfu+(\bfu-\bfgamma)\cdot\nabla\bfU,\bfH\right)\\
\hspace*{4cm}-2(\mathbb D(\bfu),\mathbb D(\bfH))-\varpi^{-1}\bfdelta\cdot\mathbb A\cdot\bfdelta+(\bff,\bfH)+  \varpi^{-1}  \bfF\cdot\bfdelta\,.
\ea
\eeq{4.46}
If we employ \eqref{1.8}--\eqref{1.10},  the property of the mollifier, assertion (v) of \lemmref{1.1}, \eqref{deH}, and recalling that $\sup_t|\bfV(t)|=1$, we can estimate the right-hand side in the following way
\be\ba{rl}\medskip
|(\bfu,\partial_t\bfH)| & \le C_2\|\nabla\bfu\|_2|\bfgamma|\,, \\ \medskip
2|(\mathbb D(\bfu),\mathbb D(\bfH))|+\lambda\left|(\bfu_\eta-{\bfgamma}+\bfU-\bfV)\cdot\nabla\bfu+(\bfu-\bfgamma)\cdot\nabla\bfU,\bfH)\right| & \le C_2|\bfdelta|\|\nabla\bfu\|_2(1+\|\nabla\bfu\|_2),\\ \medskip
|(\bff,\bfH)+  \varpi^{-1} \bfF\cdot\bfdelta| \le C_3|\bfdelta|(\|\bff\|_2+|\bfF|)& \le C_4|\bfdelta|(|\bfV|+|\dot{\bfV}|)\,.
\ea
\eeq{4.47}
We also make use of the estimate 
$$
|(\bfH,\partial_t\bfH)| \le C_2|\bfdelta||\bfgamma|.
$$
From the latter, \eqref{4.45}--\eqref{4.47} and Cauchy-Schwarz inequality we deduce that there is $\zeta_2\in(0,\zeta_1)$ such that if \be
\zeta\le\zeta_2\,, 
\eeq{Zeta}
then 
\be
\ode{E_\zeta}t+ \mbox{$\frac18$}\|\nabla\bfu\|_2^2+C_6|\bfgamma|^2 +\half\zeta\varpi^{-1}\bfdelta\cdot\mathbb A\cdot\bfdelta\le C_7\left(\zeta|\bfdelta|\|\nabla\bfu\|_2^2+|\bfV|^2+|\dot{\bfV}|^2\right)\,,
\eeq{4.48}
and, in addition, \eqref{4.44} holds. 

We claim there exists $\zeta_0\in(0,\zeta_2)$ and $\rho_0>0$ such that $E_{\zeta}(0)\le \rho_0$ implies $E_{\zeta}(T)\le \rho_0$ for $\zeta\le \zeta_0$.
%
Imposing  $E_{\zeta}(0)\le \rho_0$ with
\be
\rho_0\ge \mathcal V,
\eeq{H1}
we deduce from \eqref{ESTi} and \eqref{4.44} that
$$E(t)\le E(0)+C_2\mathcal V\le 2E_\zeta(0)+C_2\mathcal V\le (2+C_2)\rho_0,$$
yielding the estimate 
$$
|\delta(t)|\le 
C_8\rho_0^\frac12\,.
$$
Plugging this bound in \eqref{4.48}, we infer 
\be
\ode{E_\zeta}t+ \mbox{$\frac18$}\|\nabla\bfu\|_2^2+C_6|\bfgamma|^2 +\half\zeta\varpi^{-1}\bfdelta\cdot\mathbb A\cdot\bfdelta\le \frac1{16C_9}\zeta\rho_0^\frac12\|\nabla\bfu\|_2^2+C_7\left(|\bfV|^2+|\dot{\bfV}|^2\right)\,,
\eeq{4.48-bis}
and therefore
\be
\ode{E_{\zeta_0}}t+ \mbox{$\frac1{16}$}\|\nabla\bfu\|_2^2+C_6|\bfgamma|^2 +\half\zeta_0\varpi^{-1}\bfdelta\cdot\mathbb A\cdot\bfdelta\le C_7\,(|\bfV|^2+|\dot{\bfV}|^2)\,,
\eeq{4.49}
with
\be
\zeta_0:=C_9\rho_0^{-\frac12}\,.
\eeq{Zeta1}
It is classical to show, from \eqref{1.9} and H\"older inequality, that
$\|\nabla\bfu\|_2\ge C_{10}R^{-1}\|\bfu\|_2$. Therefore,  
we infer, with the help of \eqref{4.44}, 
\[
\begin{split}
\frac1{16}\|\nabla\bfu\|_2^2+C_6|\bfgamma|^2 +\half\zeta_0\varpi^{-1}\bfdelta\cdot\mathbb A\cdot\bfdelta
& \ge \frac{\zeta_0}{2}\left(C_{11}^{-1}R^{-2}\rho_0^{\frac12}\|\bfu\|_2^2 + \frac1{\varpi}\left( C_{11}^{-1}\rho_0^{\frac12}|\bfgamma|^2 +\bfdelta\cdot\mathbb A\cdot\bfdelta \right)\right)\\ 
& \ge \zeta_0 E
\ge  \half\zeta_0 E_{\zeta_0}\,,
 \end{split}
\]
if
$  \rho_0\ge C_{11}^2(R^4+1)$.
Replacing the latter into \eqref{4.49} entails
$$
\ode{E_{\zeta_0}}t+ \half\zeta_0 E_{\zeta_0}\le C_7\,(|\bfV|^2+|\dot{\bfV}|^2)\,,
 \ \ \mbox{all $t\in[\sigma,T]$}\,,
$$
which, once integrated between $\sigma$ and $T$, furnishes
$$
E_{\zeta_0}(T)\le E_{\zeta_0}(\sigma)\exp(-\half\zeta_0 T)+C_7\mathcal V\,.
$$
However, from the assertion (iii) of \lemmref{3.3} we easily show that
\be
\lim_{\sigma\to0}E_{\zeta_0}(\sigma)=E
_{\zeta_0}(0)\,
\eeq{limit}
and the claimed result that $E_{\zeta_0}(T)\le\rho_0$ follows by imposing that $\rho_0$ satisfies the condition
\be
\rho_0\ge \frac{C_7\mathcal V}{1-\exp(-\half\zeta_0 T)}\,.
\eeq{H4}
It is readily checked that \eqref{H4} certainly holds, provided we choose $\rho_0$ greater than some quantity depending only on $C_{9},C_7$, $\mathcal V$ and $T$. In fact, taking into account the choice of $\zeta_0$ made in \eqref{Zeta1} and  setting $x=C_{9}T/2\rho_0^\frac12$, \eqref{H4} is equivalent to
$$
\frac1{C_7\mathcal V}\left(\frac2{C_{9}T}\right)^2\ge x\,\frac{x}{1-{\rm e}^{-x}}\,,
$$
which is true, provided $x$ is less than a suitable quantity with the properties stated above.

Finally, the estimate \eqref{Ez} follows by integrating \eqref{4.49} over $[\sigma,T]$, letting $\sigma\to0$ and using \eqref{limit}.
\par\hfill$\square$\par

\subsection{Approximated Solutions in Bounded Domains}\label{bounded}

With the help of what we have shown so far, we are now in a position to prove the existence of a $T$-periodic weak solution to \eqref{3.1}. To this end,  we begin to give the definition of $T$-periodic weak solution to problem \eqref{3.1}. This is done exactly as we did in the case of \defref{5.1}, by replacing  $\real^3$ and $\Omega$ by $B_R$ and $\Omega_R$, respectively.
\Bd The triple $(\bfu,\bfgamma,\bfdelta)$ is a $T$-\emph{periodic weak solution} to \eqref{3.1} if
\begin{itemize} 
\item [{\rm (i)}] $\bfu \in L^2(0,T;{\mathcal D}^{1,2}(B_R))$\,, with $\bfu(x,t)|_{\partial\Omega}=\bfgamma(t)$\,, a.a. $t\in [0,T]$\,, 
$\bfgamma\in L^2(0,T;\mathbb{R}^3)$\,; 
\item [{\rm (ii)}] $\bfdelta \in W^{1,2}(0,T;\real^3)$\,;
\item [{\rm (iii)}] $(\bfu,\bfgamma,\bfdelta)$ satisfies the following equations (with $(\cdot,\cdot)\equiv(\cdot,\cdot)_{\Omega_R}$ and $\langle\cdot,\cdot\rangle\equiv\langle\cdot,\cdot\rangle_{B_R}$)
\end{itemize}
\be
\begin{split}
\Int0T\Big[-\langle\bfu,\partial_t\bfphi\rangle+\lambda\big((\bfu-{\bfgamma}+\bfU-\bfV)\cdot\nabla\bfu & + (\bfu-\bfgamma)\cdot\nabla\bfU,\bfphi)+2(\mathbb D(\bfu),\mathbb D ({\bfphi})\big)\\  & +\varpi^{-1}\hat\bfphi\cdot\mathbb A\cdot\bfdelta -(\bff,\bfphi)-\bfF\cdot \hat\bfphi\Big]{\rm d}t=0\,,\\ 
& \dot{\bfdelta}=\bfgamma\,,\ \ \Int0T\bfgamma(t){\rm d}t=\0\,,
\end{split}
\eeq{4.52}
whatever $\bfphi\in\calc_\sharp(B_R)$.
\EDD{3.1}
\Br It is easy to show that if $(\bfu,\bfgamma,\bfdelta)$ is a sufficiently regular $T$-periodic solution to \eqref{3.1}, then it satisfies \eqref{4.52} and that, with the help of the decomposition \eqref{Helm1}, the converse is also true.  
\ER{4.3}

\Bp
Let $\bfV\in W^{1,2}(0,T;\real^3)$. Then, for any $R>3R_*$ there is  at least one $T$-periodic weak solution to \eqref{3.1} in $B_R$. This solution satisfies
\be
\int_0^T(\|\nabla\bfu(t)\|_{2,\Omega_R}^2+|\bfgamma(t)|^2){\rm d}t\le C\,\mathcal V\,,
\eeq{EsTI}
where the constant $C$ is independent of $R$. Moreover, given $R_0>3R_*$, there exists a constant $C_1$ depending on $R_0$ but independent of $R$  such that for all $R\ge R_0$
\be
\int_0^T\|\bfu(t)\|_{2,\Omega_{R_0}}^2\le C_1\mathcal V\,;\ \ \ \left\|\ode{\bfu}{t}\right\|_{L^1(0,T;\cald^{-1,2}_0(\Omega_{R_0}))}\le C_1\,(\mathcal V^{\frac12}+\mathcal V)\,.
\eeq{EsTI_1} 
\EP{4.1}
{\em Proof.} Set
$$
\mathscr S_{\rho_0}:=\{(\bfu,\bfdelta)\in \calh(B_R)\times\real^3:\, \ E_{\zeta_0}(\bfu,\hat{\bfu},\bfdelta)\le \rho_0^2\}\,,
$$%
with $\zeta_0$ and $\rho_0$ as in \lemmref{4.3}. Clearly, $\mathscr S_{\rho_0}$ is a closed    convex   subset of $\calh(B_R)\times\real^3$. Next, let ${\sf s}(t):=(\bfu(t),\bfgamma(t),\bfdelta(t))
$  {be} the solution to \eqref{3.3} determined in \lemmref{3.3}, and consider the map
$$
{\sf M}: {\sf s}(0)\mapsto {\sf s}(T)\,.
$$
By \lemmref{4.3}, ${\sf M}$ maps $\mathscr S_{\rho_0}$ into itself. By \lemmref{3.3}, ${\sf M}$ is also continuous, with  $\bfu(T)\in W^{1,2}(\Omega_R)$,  thus furnishing that ${\sf M}(\mathscr S_{\rho_0})$ is compact. As a result, by Schauder fixed-point theorem we conclude that the ``mollified" problem \eqref{3.3} has at least one $T$-periodic solution.  

Our next goal is to prove that, if we let $\eta\to 0$ along a sequence $\{\eta_n\}$, the sequence of corresponding $T$-periodic (strong) solutions $(\bfu_n,\bfgamma_n,\bfdelta_n)$ converge to a $T$-periodic weak solution to \eqref{3.1}.
To this end, we remark that from \eqref{ESTi} and \eqref{Ez} it follows that
\be\ba{ll}\medskip
\Int0T(\|\nabla\bfu_n(t)\|_2^2+|\bfgamma_n(t)|^2)\,{\rm d}t\le C_1\,\mathcal V\,,\\
\Int0T\bfdelta_n(t)\cdot\mathbb A\cdot\bfdelta_n(t)\,{\rm d}t\le C_2\,\mathcal V\,,
\ea
\eeq{4.56}
where $C_1$ is independent of $R$ and $\eta$, and $C_2$  depends on $R$ but is independent of $\eta$.
   Multiplying both sides of \eqref{3.3}$_1$ --written for these solutions-- by the test function $r(t)\bfpsi$, we show that for all $\bfpsi\in\cald^{1,2}(B_R)$ and all smooth $r$, 
\be
\begin{split}
\langle\bfu_n(t),r(t)\bfpsi\rangle-\langle\bfu_n(0),r(0)\bfpsi\rangle=-\Int0t\Big[-\langle\bfu_n,r^\prime\bfpsi\rangle+\lambda\left(((\bfu_n)_{\eta_n}-{\bfgamma_n}+\bfU-\bfV)\cdot\nabla\bfu_n\right.\\ \medskip\left.+(\bfu_n-\bfgamma_n)\cdot\nabla\bfU,r\bfpsi\right)+\left.2(\mathbb D(\bfu_n),\mathbb D ({r\bfpsi}))+\varpi^{-1}r{\hat{\bfpsi}}\cdot\mathbb A\cdot\bfdelta_n-(\bff,r\bfpsi)-\bfF\cdot r\hat{\bfpsi}\right]{\rm d}s\,.
\end{split}
\eeq{4.56_1}
Now, by a standard procedure (see e.g. \cite[Section 2]{Gaintro}), we use this relation for $t=T$, \lemmref{2.3} and the functional properties of $(\bfu_n,\bfgamma_n,\bfdelta_n)$ to deduce that for all $\bfphi\in\calc_\sharp(B_R)$, it holds
\be
\ba{ll}\medskip
\Int0T\big[-\langle\bfu_n,\partial_t\bfphi\rangle+\lambda\left(((\bfu_n)_{\eta_n}-{\bfgamma_n}+\bfU-\bfV)\cdot\nabla\bfu_n+(\bfu_n-\bfgamma_n)\cdot\nabla\bfU,\bfphi\right)\\ \ms\hspace*{2.5cm}+\left.2(\mathbb D(\bfu_n),\mathbb D ({\bfphi}))+\varpi^{-1}\hat\bfphi\cdot\mathbb A\cdot\bfdelta_n-(\bff,\bfphi)-\bfF\cdot \hat\bfphi\right]{\rm d}s=0\,.
\ea
\eeq{4.57} 
We also recall from \eqref{wefo}$_2$:
\be
\dot{\bfdelta}_n(t)=\bfgamma_n(t)\,,\ \ t\in[0,T]\,,\,\ \ \Int0T\bfgamma_n(t){\rm d}t=\0.
\eeq{4.58}
By \eqref{4.56}, \eqref{4.58} and \eqref{1.10}, we infer that there exist
$$
(\bfu,\bfgamma,\bfdelta)\in L^2(0,T;\cald^{1,2}(B_R))\times L^{2}(0,T;\real^3)\times W^{1,2}(0,T;\real^3)\,,
$$ 
with $\hat{\bfu}(t)=\bfgamma(t)$ such that, as $n\to\infty$,\footnote{In what follows we will not make  notational distinction between sequences and subsequences.}
\be\ba{ll}\medskip
\bfu_n\to\bfu\,,\ \ \mbox{weakly in $L^2(0,T;\cald^{1,2}(B_R))$}\,,
\\ \medskip
\bfgamma_n\to\bfgamma\,,\ \ \mbox{weakly in $L^2(0,T;\real^3)$},\\
\bfdelta_n\to \bfdelta\,,\ \ \mbox{weakly in $W^{1,2}(0,T;\real^3)$,\,  and in $C([0,T];\real^3)$\,.}
\ea
\eeq{4.59}
From \lemmref{4.4}, \eqref{4.56}, \eqref{4.59}$_1$ and Simon compactness theorem \cite{Simo}, we also get
\be
\bfu_n\to\bfu\,,\ \ \mbox{strongly in $L^2(0,T;\calh(B_R))$}\,,
\eeq{4.60}
which  implies
\be
\bfgamma_n\to\bfgamma\,,\ \ \mbox{strongly in $L^2(0,T;\real^3)$}\,.
\eeq{4.60_2}
We now pass to the limit $n\to\infty$ in \eqref{4.57}--\eqref{4.58}. Employing \eqref{4.59}--\eqref{4.60_2}, it is not difficult to show that, in doing so, we can replace everywhere in \eqref{4.57}--\eqref{4.58} $\bfu_n,\bfgamma_n$ and $\bfdelta_n$ with $\bfu,\bfgamma$ and $\bfdelta$, respectively, with \eqref{4.58}$_1$ holding for a.a. $t\in[0,T]$. The only point that deserves a little care is the convergence of the nonlinear term:
\be
I_n:=\int_0^T(((\bfu_n)_{\eta_n}-{\bfgamma_n})\cdot\nabla\bfu_n,\bfphi){\rm d}t\to \int_0^T((\bfu-\bfgamma)\cdot\nabla\bfu,\bfphi){\rm d}t=:I\,.
\eeq{4.61}
To show \eqref{4.61}, we first observe that, using Schwarz inequality, 
$$\ba{rl}\medskip
|I_n-I|\le &\!\!\!\Max{t,x}|\bfphi(x,t)|\Int0T\left(\|(\bfu_n-\bfu)_{\eta_n}\|_2^2+\|((\bfu)_{\eta_n}-\bfu\|_2^2+ |\bfgamma_n-\bfgamma|^2\right)\|\nabla\bfu_n\|_2^2{\rm d}t\\
&\!\!\!+\left|\Int0T((\bfu-\bfgamma)\cdot(\nabla\bfu_n-\nabla\bfu),\bfphi){\rm d}t\right|\,.
\ea
$$
Then, using \eqref{4.59}$_1$, \eqref{4.60}, \eqref{4.60_2} and classical properties of the mollifier,  we infer the convergence in \eqref{4.61} holds. Finally, \eqref{EsTI} is established by letting $n\to\infty$ in \eqref{4.56}$_1$ and using \eqref{4.59}$_{1,2}$. 

It remains to prove \eqref{EsTI_1}. The first inequality in \eqref{EsTI_1} is an obvious consequence of \eqref{1.9} and \eqref{EsTI}. To show the second one,  we choose in  \eqref{4.52} $\bfphi=\psi\bfhi$, for arbitrary $\bfhi\in\calc_0(B_{R_0})$ and $\psi\in C_0^\infty((0,T);\real)$. We thus obtain 
\be
\int_0^T\psi^\prime(t)(\bfu(t),\bfhi)\,{\rm d }t=\int_0^TG_{\mbox{\tiny $\bfhi$}}(t)\psi(t)\,{\rm d}t\,,
\eeq{4.66}
where
$$
G_{\mbox{\tiny $\bfhi$}}(t)=-\lambda\left((\bfu-{\bfgamma}+\bfU-\bfV)\cdot\nabla\bfu+(\bfu-\bfgamma)\cdot\nabla\bfU,\bfhi\right)-2(\mathbb D(\bfu),\mathbb D ({\bfhi}))+(\bff,\bfhi)
\,.
$$
Arguing exactly as in the proof of \lemmref{4.4}, we then prove that there is a constant $C$ depending on $R_0$ such that
\be
|G_{\mbox{\tiny $\bfhi$}}(t)|\le  C(R_0)\, (|\bfV(t)|+\|\nabla\bfu(t)\|_2+\|\nabla\bfu(t)\|_2^2)\|\mathbb D(\bfhi)\|_2\,,\ \ \mbox{for all $\bfhi\in \calc_0(B_{R_0})$}.
\eeq{4.67}
Henceforth, since $\calc_0(B_{R_0})$ is dense in $\cald_0^{1,2}(\Omega_{R_0})$, we infer
$$
G_{\mbox{\tiny $\bfhi$}}(t)=[\bfg(t),\bfhi]
$$
for some $\bfg(t)\in \cald_0^{-1,2}(\Omega_{R_0})$, where $[\cdot,\cdot]$     denotes the duality pairing $\cald_0^{-1,2}(\Omega_{R_0})\leftrightarrow \cald_0^{1,2}(\Omega_{R_0})$. 
This, in combination with \eqref{4.66} and the arbitrariness of $\psi$, furnishes
$$
\ode{}t(\bfu,\bfhi)=[\bfg(t),\bfhi]\,,
$$
in the sense of distribution.
Furthermore, in view of \eqref{EsTI} and \eqref{4.67}, we have
$$
\int_0^T|G_{\mbox{\tiny $\bfhi$}}(t)|{\rm d}t\le C\left(\mathcal V^{\frac12}+\mathcal V\right)\|\mathbb D(\bfhi)\|_2\,,
$$  
and the desired property is then proved.
\par\hfill$\square$\par
\subsection{$T$-periodic Weak Solutions for the Original Problem}\label{sec:Fin} 
This last subsection is dedicated to the proof of \theoref{5.1}.
Let $\{\Omega_n\equiv\Omega_{R_n}\}$, $R_1>3R_*$, be a sequence of ``invading domains," namely,
$$
\Omega_{n-1}\subset\Omega_{n}\,,\ n\in\nat\,;\ \,\ \cup_{n=1}^\infty\Omega_n=\Omega\,,
$$
and let $\{{\sf s}_n\equiv (\bfu_n,\bfgamma_n,\bfdelta_n)\}$ be the sequence of corresponding $T$-periodic weak solutions determined in \propref{4.1}. For each $n$, we extend $\bfu_n$ to 0 outside $\Omega_n$ and continue to denote by $\bfu_n$ its extension. Consequently, by \remref{2.1}, $\{\bfu_n\}\subset W^{1,2}(\Omega)\cap \cald^{1,2}(\real^3)$. From the bound \eqref{EsTI}, we deduce that there is a subsequence of $\{(\bfu_n,\bfgamma_n)\}$,
again denoted by the same symbol, with $\bfu_n|_{\Omega_0}=\bfgamma_n$, and functions $(\bfu,\bfgamma)\in L^2(0,T; \cald^{1,2}(\real^3))\times L^2(0,T;\real^3)$ such that 
\be\ba{ll}\medskip
\bfu_n\to\bfu\,,\ \ \mbox{weakly in $L^2(0,T;\cald^{1,2}(\real^3))$}\,,
\\ 
\bfgamma_n\to\bfgamma\,,\ \ \mbox{weakly in $L^2(0,T;\real^3)$}\,,
\ea
\eeq{5.3} 
and for which \eqref{5.2} holds.
Fix $R_0 > R_*$ arbitrarily. From  \eqref{EsTI_1}, \eqref{5.3}$_1$ and Aubin-Lions-Simon theorem \cite{Simo},  we can extract another subsequence, again denoted by $\{\bfu_n\}$, such that
\be
\bfu_n\to\bfu\,,\ \ \mbox{strongly in $L^2(0,T;L^{2}(\Omega_{R_0}))$}\,.
\eeq{5.4}
Also, from a classical trace inequality, we get
$$
|\bfgamma_n-\bfgamma|\le c_\epsilon\,\|\bfu_n-\bfu\|_2+\epsilon\,\|\nabla(\bfu_n-\bfu)\|_2\,,
$$
where $\epsilon>0$ is arbitrary \cite[Exercise II.4.1]{Gab}, so that, from the latter, \eqref{5.3} and \eqref{5.4} we deduce
\be
\bfgamma_n\to\bfgamma\,,\ \ \mbox{strongly in $L^2(0,T;\real^3)$}\,.
\eeq{5.5}
Since for all $n\in\mathbb N$ we have
$$\Int0T\bfgamma_n(t){\rm d}t=\0\,,$$
we deduce that $\bfgamma$ has zero average.
Now set 
$$
\tilde{\bfdelta}_n:=\bfdelta_n-\bar{\bfdelta_n}\,;\ \ \ \Big(\bar{w}:=T^{-1}\int_0^Tw(t)\,{\rm d}t\Big) \,.
$$ 
Recalling that $\dot{\bfdelta_n}=\bfgamma_n$, Poincar\'e-Wirtinger  inequality yields
$$
\int_0^T|\tilde\bfdelta_n-\tilde\bfdelta_m|^2\le T^2\int_0^T|\bfgamma_n-\bfgamma_m|^2\,,
$$
so that $\{\tilde\bfdelta_n\}$ is a Cauchy sequence in $L^2(0,T;\real^3)$.
 It follows that there exists $\tilde\bfdelta\in W^{1,2}(0,T;\real^3)$ with zero average such that $\dot{\tilde\bfdelta}=\bfgamma$ and
\be
\tilde{\bfdelta}_n\to \tilde{\bfdelta} \,, \quad \mbox{strongly in }W^{1,2}(0,T;\real^3).
\eeq{5.6}
The last subsequence we have selected  may depend on $R_0$.
However, covering $\Omega$, with an increasing sequence of bounded domains and using
Cantor diagonal method, we may extract a further subsequence  for which all the
above properties, and in particular \eqref{5.4}, hold for {\em all} $R_0$. 

In order to complete the proof of the theorem, it remains to show that the limiting functions determined above satisfy 
 the weak formulation of \eqref{3.1}. In particular, we still need to prove the convergence of the sequence of averages $\{\bar\bfdelta_n\}$. 
From the weak formulation \eqref{4.52}$_1$ satisfied by  $(\bfu_n,\bfgamma_n,\bfdelta_n)$ and the arbitrariness of $R$ it follows that \mbox{for any  fixed at will $\bfphi\in\calc_\sharp(\real^3)$} the sequence $\{{\sf s}_n\}$ obeys the following equation for all sufficiently large $n\in\nat$
\be\ba{ll}\medskip
\Int0T\big[\langle\bfu_n,\partial_t\bfphi\rangle-\lambda\left((\bfu_n-{\bfgamma}_n+\bfU-\bfV)\cdot\nabla\bfu_n+(\bfu_n-\bfgamma_n)\cdot\nabla\bfU,\bfphi\right)\\ \hspace*{2.5cm}\left.-2(\mathbb D(\bfu_n),\mathbb D ({\bfphi}))-\varpi^{-1}\hat\bfphi\cdot\mathbb A\cdot\tilde{\bfdelta}_n+(\bff,\bfphi)+\bfF\cdot \hat\bfphi\right]{\rm d}t=\varpi^{-1}T\,\bar{\hat\bfphi}\cdot\mathbb A\cdot\bar{\bfdelta_n}\,.
\ea
\eeq{5.7}
The convergences proved for the sequences $\{\bfu_n,\bfgamma_n,\tilde\bfdelta_n\}$ in \eqref{5.3}$_1$, \eqref{5.4} --valid for {\em all} $R_0$-- and  \eqref{5.5}--\eqref{5.6}, and an argument similar to that used in the proof of \eqref{4.61}, we show that, as $n\to\infty$, the left-hand side of \eqref{5.7} converges to the same quantity with $(\bfu_n,\bfgamma_n,\tilde\bfdelta_n)$ replaced by $(\bfu,\bfgamma,\tilde\bfdelta)$, for any arbitrarily fixed $\bfphi$. This implies that $\bar{\bfdelta_n}\to\bar{\bfdelta}$, for some $\bar{\bfdelta}\in\real^3$. 
We then deduce that the triple $(\bfu,\bfgamma,\bfdelta=\tilde\bfdelta+\bar{\bfdelta})$ satisfies \eqref{4.521}$_1$ whatever the test function $\bfphi$ taken in $\calc_\sharp(\real^3)$.
This concludes the proof. 
\par\hfill$\square$\ms\par\noindent

\section*{Declarations}
\begin{itemize}
\item The authors have no relevant financial or non-financial interests to disclose.
\item The authors have no competing interests to declare that are relevant to the content of this article.
\item Data sharing not applicable to this article as no datasets were generated or analysed during the current study.
\end{itemize}

\section*{Acknowledgments} 
This work is  partially supported by the National Science Foundation Grant DMS--2307811 and the ARC Advanced 2020-25 ``PDEs in interaction''. D. Bonheure is supported by the Francqui Foundation as Francqui Research Professor 2021-24.

\ed

\ed